\documentclass[11pt]{amsart}
\usepackage{amssymb,amsmath,amsfonts,amscd,color,euscript,graphics,comment}
\usepackage{enumerate}
\usepackage{amssymb}
\usepackage{mathtools}
\usepackage{amsthm}
\usepackage{latexsym}
\usepackage{multicol}
\usepackage{verbatim}
\usepackage[all]{xy}
\usepackage{hyperref}

\advance\textwidth by 1.2in \advance\oddsidemargin by -.6in \advance\evensidemargin by -.6in

\newcommand{\eval}{\text{ev}}
\newcommand{\supp}{\text{supp}}

\newcommand{\Id}{\text{Id}}
\newcommand{\bc}{\mathbb{C}}
\newcommand{\bz}{\mathbb{Z}}
\newcommand{\bu}{\mathbf U}
\newcommand{\wt}{\operatorname{wt}}
\newcommand{\Ext}{\operatorname{Ext}}
\newcommand{\Hom}{\operatorname{Hom}}
\newcommand{\Ob}{\operatorname{Ob}}
\newcommand{\ann}{\operatorname{Ann}}
\newcommand{\id}{\operatorname{Id}}

\newcommand{\mode}{\operatorname{mod}}

\newcommand{\sym}{\operatorname{sym}}
\newcommand{\rank}{\operatorname{rank}}
\newcommand{\adm}{\operatorname{adm}}
\newcommand{\br}{\mathbf R}
\newcommand{\bw}{\mathbf W}
\newcommand{\ba}{\mathbf A}
\newcommand{\bop}{\mathbf p}
\newcommand{\Max}{\operatorname{Max}}
\newcommand{\Spec}{\operatorname{Spec}}
\theoremstyle{definition}

\theoremstyle{definition}
\newtheorem*{defn}{Definition}
\theoremstyle{definition}
\newtheorem*{thm}{Theorem}
\newtheorem*{theom}{Theorem}

\newtheorem*{rem}{Remark}

\newtheorem*{lem}{Lemma}
\newtheorem*{cor}{Corollary}
\newtheorem*{prop}{Proposition}

\newcounter{cnt}

\newenvironment{pf}{\proof}{\endproof}
\newenvironment{enumerit}{\begin{list}{{\hfill\rm(\roman{cnt})\hfill}}{%
\settowidth{\labelwidth}{{\rm(iv)}}\leftmargin=\labelwidth%
\advance\leftmargin by \labelsep\rightmargin=0pt\usecounter{cnt}}}{\end{list}} \makeatletter
\def\mydggeometry{\makeatletter\dg@YGRID=1\dg@XGRID=20\unitlength=0.003pt\makeatother}
\makeatother \theoremstyle{remark}

\numberwithin{equation}{section}

\begin{document}

\renewcommand{\frak}{\mathfrak}
\newcommand{\nc}{\newcommand}
\newcommand{\rnc}{\renewcommand}
\newcommand{\lie}[1]{\mathfrak{#1}}
\newcommand{\liefixed}[1]{\mathfrak{#1}_0}
\newcommand{\laurent}{\bc\left[t^{\pm 1}\right]}

\newcommand{\bm}{\mathbf{m}}
\newcommand{\Lgg}{L^\Gamma(\lie g)}
\newcommand{\Lhg}{L^\Gamma(\lie h)}
\newcommand{\Lat}{L^{\Gamma}(\lie{sl}_3)}
\newcommand{\Lg}{L(\lie g)}
\newcommand{\Lh}{L(\lie h)}
\newcommand{\Bagg}{\ba_{\lambda}^{\Gamma}}
\newcommand{\cal}{\mathcal}
\newcommand{\Maps}{Maps}
\newcommand{\ev}{ev}

\title
{Global Weyl modules for the twisted loop algebra}
\author{Ghislain Fourier, Nathan Manning and Prasad Senesi}
\address{Ghislain Fourier:\newline
Mathematisches Institut, Universit\"at zu K\"oln, Germany}
\email{gfourier@math.uni-koeln.de}
\address{Nathan Manning:\newline
University of California, Riverside}
\email{nmanning@math.ucr.edu}
\address{Prasad Senesi:\newline
The Catholic University of America}
\email{senesi@cua.edu}
\begin{abstract}
We define global Weyl modules for twisted loop algebras and analyze their highest weight spaces, which are in fact isomorphic to Laurent polynomial rings in finitely many variables. We are able to show that the global Weyl module is a free module of finite rank over these rings. Furthermore we prove, that there exist injective maps from the global Weyl modules for twisted loop algebras into a direct sum of global Weyl modules for untwisted loop algebras. Relations between local Weyl modules for twisted and untwisted generalized current algebras are known; we provide for the first time a relation on global Weyl modules.
\end{abstract}

\maketitle \thispagestyle{empty}
\section{Introduction}

Let $\lie g$ be a simple compex Lie algebra.  The representation theory of the loop algebras $\lie g~\otimes~\bc[t^{\pm1}]$ (denoted by $L(\lie g)$) has been an active area of research for several decades.  The global Weyl modules, initially defined in \cite{CP01}, have played a prominent role in this theory since their introduction over ten years ago.  

The global Weyl modules are naturally indexed by the dominant integral weights $P^+$ of $\lie g$.  For such a weight $\lambda$, the corresponding global Weyl module $W(\lambda)$ was originally defined as a certain maximal integrable highest weight module generated by a non-zero vector of weight $\lambda$. But this module can also be defined as a projective object in a category of locally finite representations for $\Lg$ (\cite{CFK10}).   
 
The study of Weyl modules was motivated by the representation theory of quantum affine algebras, and in fact certain quotients of the global Weyl modules (the \textit{local Weyl modules}) are $q=1$ limits of simple representations of quantum affine algebras.  But these local Weyl modules are not simple as modules for the loop algebra in general.  In fact, the category of finite--dimensional modules for $\Lg$ is not semi--simple, and the local Weyl modules are in some sense the largest indecomposable highest weight modules. The calculation of their dimension and character has led to a series of papers (\cite{CP01}, \cite{CL06}, \cite{FoL07}, \cite{Na11}, \cite{BN04}). 

Although the global Weyl module $W(\lambda)$ is an infinite--dimensional module, it was shown in \cite{CP01} that its highest weight space $\ba_\lambda$ is isomorphic to a Laurent polynomial ring in finitely many variables, and any local Weyl module can be obtained from a global Weyl module by tensoring the global Weyl module with $\ba_\lambda / \mathbf{m}$, where $\bm$ is a maximal ideal of $\ba_\lambda$.  This construction justifies the terminology \textit{local} and \textit{global} Weyl module.  In the aforementioned series of papers, it was shown that the dimension and $\lie g$-character of a local Weyl module is independent of the maximal ideal; hence the global Weyl module is a free module of finite rank for the algebra $\ba_{\lambda}$.

We describe several ways to generalize the notion of a global Weyl module, all obtained by generalizing the Lie algebras $\lie g \otimes \bc[t^{\pm 1}]$.  The algebra $\bc[t^{\pm1}]$ can be replaced with a more general algebra. The first step in this direction was made in \cite{FL04}, where global Weyl modules for $\lie g \otimes A$ are studied when $A$ is the coordinate ring of an affine variety (the \textit{generalized current algebra} $\lie g \otimes A$ can be viewed here as the Lie algebra of maps from $\Spec A$ to $\lie g$).  It was shown here that the global Weyl modules exist and are non-zero. In \cite{CFK10}, this direction was further generalized, by taking for $A$ any commutative, associative and unital algebra over $\bc$.  Global Weyl modules and the algebras $\ba_\lambda$ were defined, and  homological methods were developed for the further analysis of all of these modules. 
In this paper (and also in \cite{FL04}) it was shown that the global Weyl module is in general not a free right $\ba_\lambda$--module - this may be a special property of the case when $A=\bc[t^{\pm 1}]$.

We mention several other generalizations of the loop algebras for which a global Weyl module could be studied.  
The \textit{twisted loop} (\textit{sub})\textit{algebras} of $\lie g \otimes \bc[t^{\pm 1}]$, here denoted $\Lgg$, are the algebras of fixed points in $\lie g \otimes \bc[t^{\pm1}]$ under a group action of $\Gamma = \bz/m \bz$ obtained from an order $m$ Dynkin diagram automorphism of $\lie g$ ($\Gamma$ acts upon $\bc^*$ by multiplication with an $m^{th}$ primitive root of unity; using the induced action on the coordinate ring, we have a diagonal action on $\lie g \otimes \bc[t^{\pm1}]$). 
In \cite{CFS08}, local Weyl modules were defined and studied for these algebras. The main result was that any local Weyl module of $\Lgg$ can be obtained by restricting a local Weyl module for $\Lg$. 

More generally, let $\Gamma$ be any finite group, acting on $\lie g$ and an affine scheme $X$ by automorphisms.  The Lie algebras $(\lie g \otimes A)^{\Gamma}$ of equivariant regular maps from $X$ to $\lie g$ are called the equivariant map algebras. Their finite--dimensional simple modules were studied and classified in \cite{NSS}.  Local Weyl modules for an equivariant map algebra were defined and studied in \cite{FKKS11} in the case when $\Gamma$ is abelian and its action on $X$ is free.  Analogous to the case of twisted loop algebras, it was shown again there that any local Weyl module for $(\lie g \otimes A)^{\Gamma}$ can be obtained by restriction from a local Weyl module for $\lie g \otimes A$.

Let $A = \bc[t]$ and $\Gamma = \bz /m \bz$ as above, acting upon $\bc$ by a primitive $m^{th}$ root of unity.  When $\Gamma \neq 1$, this action is not free, and the results of \cite{FKKS11} do not apply.  Local Weyl modules were studied and defined in \cite{FK11}. It was shown that there exist local Weyl modules (supported at the origin) which have no counterpart in the category of modules for $\lie g \otimes \bc[t]$.

All of the above mentioned papers concerning fixed--point algebras have studied only finite--dimensional representations and, in particular, only the local Weyl modules.  The notion of a global Weyl module for all of these algebras has never before been approached.

In this paper, we define for any dominant integral weight $\lambda$ of the fixed-point algebra $\lie g _0$ of $\lie g$ a global Weyl module $W^{\Gamma}(\lambda)$.  We also describe its highest weight space $\Bagg$, by giving it a natural algebra structure, and define a Weyl functor $\bw^{\Gamma}_\lambda$ from the category of left $\Bagg$-modules to the category of integrable $\Lgg$-modules. As in \cite{CFK10}, we obtain, by using the Weyl functor, a homological characterization of Weyl modules. We prove that $\Bagg$ is a finitely generated ring of symmetric Laurent polynomials. We identify its maximal spectrum with equivariant finitely supported functions from $\bc^*$ to $P^+$ (analogous to \cite{NSS}) as well as with certain finite multisets. Furthermore, we prove there is a canonical embedding of $\Bagg$ into $\ba_{\mu}$ for all $\mu \in P^+$ satisfying $\mu\left.\right|_{\lie h_0} = \lambda$, where $\lie h_0$ is the fixed--point subalgebra of a Cartan subalgebra $\lie h$ of $\lie g$. Using results from \cite{CFS08}, we show that $W^{\Gamma}(\lambda)$ is a free module for $\Bagg$.
The global Weyl module $W(\mu)$ ($\mu \in P^+$) for the loop algebra $\Lg$ is, via restriction, also a module for $\Lgg$. Furthermore, if $\lambda$ is a dominant integral weight of $\lie g_0$, then
$ \bigoplus_{\mu} W(\mu)$ (where the sum is taken over all $\mu \ : \mu\left. \right|_{\lie h_0}= \lambda$) is a $\Lgg$-module and our main theorem compares this module with $W^{\Gamma}(\lambda)$:
\begin{theom}
Let $\lambda$ be a dominant integral weight of $\lie g_0$.  Then there exists a natural embedding of $\Lgg$-modules
$$W^\Gamma(\lambda) \hookrightarrow \bigoplus_{\mu } W(\mu),$$
where the sum is over all $\mu \in P^+$ satisfying $ \mu\left. \right|_{\lie h_0}= \lambda$.
\end{theom}
This result might suggest an approach to define global Weyl modules in a more general setting. A significant obstruction to do so for equivariant map algebras is the absence (in general) of a non-zero Cartan subalgebra; hence the notion of weights is unavailable.

Some motivation for our work also comes from the finite-dimensional representation theory of the quantum affine algebra, where relationships are known (\cite[Theorem 4.15]{Her10}) between Kirillov-Reshetikhin modules for the twisted and untwisted algebras. 

This paper is organized as follows. In Section~\ref{section2} we give basic definitions and notation, and recall the known facts about simple modules for $\Lg$ and $\Lgg$. In Section~\ref{section3} we define our main object of study, the global Weyl module, and state the main result (Theorem~\ref{maintheorem}). Section~\ref{section4} is dedicated to the analysis of the Weyl functor and its properties; it is also shown that the global Weyl module is finitely generated as a module for the finitely generated algebra $\Bagg$. In Section~\ref{section-bagg}, the algebra $\Bagg$ is analyzed and identified with a ring of symmetric Laurent polynomials. Section~\ref{local-comp} recalls the various definitions of local Weyl modules and gives the proof that the global Weyl module is free as an $\ba^\Gamma_\lambda$--module.  Finally in Section~\ref{section8} the main theorem (Theorem~\ref{maintheorem}) is proven.

\vskip12pt
\textbf{Acknowledgements:} 
The first and second author would like to thank the Hausdorff Research Institute for Mathematics and the organizers of the Trimester Program on the Interaction of Representation Theory with Geometry and Combinatorics, during which many of the ideas in the current paper were developed. The authors would like to thank Alistair Savage for helpful comments, while the first author would also like to thank Deniz Kus for stimulating discussions. The first and third author are grateful to Vyjayanthi Chari for her hospitality at UC Riverside, when parts of this paper were written. The second author is grateful to David Hernandez for his hospitality and for early discussions at Paris 7. The first author was partially sponsored by the DFG-Schwerpunktprogramm 1388 ``Darstellungstheorie".


\section{Preliminaries}\label{section2}
\subsection{}
 Throughout the paper $\bc$  denotes the set of complex  numbers and
$\bz_+$ the set of non--negative integers.  Let  $\lie g$   be a finite--dimensional
   simple Lie algebra of rank $n$ with Cartan matrix $(a_{ij})_{i,j\in I}$ where $I=\{1,\cdots, n\}$. Fix a Cartan subalgebra $\lie h$ of $\lie g$  and let
  $R$ denote the corresponding  set of   roots.  Let $\{\alpha_i\}_{i\in I}$ (resp. $\{\omega_i\}_{i\in I}$)  be  a
set of simple roots (resp. fundamental weights) and  $Q$ (resp. $Q^+$), $P$ (resp. $P^+$) be the integer span (resp. $\bz_+$--span) of the
simple roots and fundamental weights respectively. Denote by  $\le $  the usual partial order on $P$,  $$\lambda,\mu\in
P,\ \ \lambda\le \mu\ \iff\  \mu-\lambda\in Q^+.$$ Set $R^+= R\cap Q^+$ and let  $\theta$ be  the unique maximal element in $R^+$ with
respect to the partial order.

 We fix a Chevalley basis  $X^\pm_\alpha$, $H_i$, $\alpha\in R^+$, $i\in I$ of $\lie g$ and set $X_i^\pm=X^\pm_{\alpha_i}$,
 $H_\alpha=[X^+_\alpha, X^-_\alpha]$ and note that $H_i=H_{\alpha_i}$. For each $\alpha\in R^+$, the  subalgebra of $\lie g$ spanned by $\{X^\pm_\alpha, H_\alpha\}$ is isomorphic to $\lie{sl}_2$.
 Define subalgebras  $\lie n^\pm $ of $\lie g$
  by  $$\lie n^\pm=\bigoplus_{\alpha\in R^+}\bc X^\pm_\alpha,$$
   and note that $$\lie g=\lie n^-\oplus\lie h\oplus\lie n^+.$$ Given any Lie algebra $\lie a$,
    let $\bu(\lie a)$ be the universal enveloping algebra of $\lie a$.  The map
    $x\mapsto x\otimes 1+1 \otimes x$, $x\in\lie a$ extends to an algebra homomorphism $\Delta: \bu(\lie a)\to \bu(\lie a)\otimes\bu(\lie a)$.
By the Poincare--Birkhoff--Witt theorem, we know that if $\lie b$ and $\lie c$ are Lie subalgebras of $\lie a$ such that   $\lie a=\lie b\oplus\lie c$ as vector spaces, then $$\bu(\lie a)\cong\bu(\lie b)\otimes\bu(\lie c)$$ as vector spaces.

Let $\sigma$ be a diagram automorphism of $\lie g$ of order $m$ and take $\Gamma=\left< \sigma \right> = \bz/m\bz$. Fix $\zeta$ a primitive $m^{th}$ root of 1. Then $\lie g$ decomposes as 
\[ \lie g=\bigoplus_{s=0}^{m-1} \lie g_{s} \]  where
\[ \lie g_{s}=\left\{x\in g\ :\  \sigma(x)=\zeta^s x\right\}.\]
This provides a $\Gamma$-grading of $\lie g$.  Given any subalgebra $\lie a$ of $\lie g$ which is preserved by $\Gamma$,  set $\lie a_{s} = \lie g_{s}\cap \lie a$.  The following is well known (see for example \cite{Ca05} or \cite[Chapter 8]{K90}). $\liefixed g$ is a simple Lie algebra and $\liefixed h$ is a Cartan subalgebra of $\liefixed g$, and we denote by $R_0$ the corresponding set of roots.  We fix a set of simple roots $\left\{ \alpha_i \right\}_{i \in I_0}$, and let $Q_0$ (resp. $Q_0^+$), $P_0$ (resp. $P_0^+$) be the integer span (resp. $\bz_+$--span) of the simple roots $\left\{ \alpha_i \right\}_{i \in I_0}$ and the weights $\{\omega_i\}_{i \in I_0}$ (resp. in the case where $\lie g$ is of type $A_{2n}$ the span of $\{\omega_i\}_{i \in I_0\setminus \{\rank \lie g_0\}} \cup \{2 \omega_{\rank \lie g_0}\}$), respectively. Moreover,  $\lie g_{s}$ is an irreducible representation of $\liefixed g$ for all $s$, and 
\[ \lie n^{\pm} \cap \liefixed g = \lie  n_0^{\pm} = \bigoplus_{\alpha \in R_0^+}(\lie g)_{\pm \alpha}.\]
We set $\lie h_0 = \lie h \cap  \lie g_0$, so we have $\lie g = \lie n_0^- \oplus \lie h_0 \oplus \lie n_0^+$ is a triangular decomposition of $\lie g_0$.  The rank of $\liefixed g$ is equal to the number of orbits of $I$ under the induced action of $\Gamma$.  We identify this set of orbits with an index set for the simple roots of $\liefixed g$, and further identify this with a subset $I_0 = \{1, \ldots, \text{rk } \liefixed g \} \subset I$ by choosing a representative of each orbit.
 
As $\lie h_0 \subseteq \lie h$, we have a natural map $\lie h^* \twoheadrightarrow \lie h_0^*$; furthermore our choice of Chevalley basis elements $\left\{ h_i \right\}$ is such that $P \twoheadrightarrow P_0$.  If $\lambda \in P$, we denote its image under this projection by $\bar{\lambda} \in P_0$.  It will sometimes become convenient to label elements in $P_0$ as images of this projection, and if not, the context will clarify whether a functional $\lambda$ lives in $P$ or $P_0$.  

As diagram automorphism, the group $\Gamma$ acts upon the nodes $\left\{ 1, \ldots, \text{rk}(\lie g) \right\}$ of the Dynkin diagram of $\lie g$, and for a node $i$ of this diagram we denote by $\Gamma_i$ the stabilizer of $i$ in $\Gamma$. More generally, $\Gamma$ acts on $R$ and we denote by $\Gamma_{\alpha}$ the stabilizier of $\alpha$. For $\alpha \in R_0$, we often choose a preimage lying in  $R$, and when this will not cause confusion, we also label it $\alpha$. For $0 \leq k < m$ and $\alpha \in R_0$, we define the following elements $h_\alpha(k) \in \lie h \cap \lie g_k$, $x^\pm_\alpha(k) \in \lie n^\pm \cap \lie g_k$:

\[ h_\alpha(k) = \frac{1}{\left| \Gamma_\alpha \right| }\sum_{j=0}^{m-1} (\zeta^k)^j H_{\sigma^j(\alpha)}, \; \; \; x_\alpha^\pm(k) = \frac{1}{\left| \Gamma_\alpha \right| } \sum_{j=0}^{m-1} (\zeta^k)^j X^\pm_{\sigma^j(\alpha)}.\]

If $\lie g$ is of type $A_{2n}$ and $\alpha \in R_0$ is a short root, we have the additional elements $x_{2\alpha}^\pm(1)$ with 
$$\bc x_{2\alpha}^\pm(1) = \bc [ X^\pm_{\alpha}, X^\pm_{\sigma(\alpha)}].$$

Observe, that if $\Gamma_\alpha = \Gamma$, then $h_\alpha(0) = h_\alpha(\epsilon)$ for all $0 \leq \epsilon < m$. We set $h_i(k) \coloneqq h_{\alpha_i}(k)$, $x_i^{\pm}(k) \coloneqq x^{\pm}_{\alpha_i}(k)$ and observe $h_i = h_i(0)$, $x_i^{\pm} = x_i^{\pm}(0)$.\\
Then  for all $x^\pm_\alpha(0) \in (\liefixed g)_{ \alpha}$, $h_\alpha(0) \in \lie h$, the subalgebra generated by 
$\left\{ x^\pm_\alpha(0), h_\alpha(0) \right\}$ generate a Lie algebra isomorphic to $\lie{sl}_2$, and $\left\{x^\pm_\alpha(0), h_i\right\}_{i \in I_0, \alpha \in R_0^+}$ is a Chevalley basis of $\liefixed g$; see \cite{FK11} for details. In the case when $\Gamma = \Id$, we have $\liefixed g = \lie g$, $x_\alpha(0) = X_\alpha$, $h_\alpha(0) = H_\alpha$, $P_0^{\pm} = P^{\pm}$ and $Q_0^{\pm} = Q^{\pm}$.

\subsection{}\label{subsectiongeometry} 
Let $A=\laurent$ and let $A_+$ be a fixed vector space complement to the subspace $\bc$ of $A$.  Given a Lie algebra $\lie a$,  define a Lie algebra structure on $\lie a\otimes A$,  by 
$$[x\otimes a, y\otimes b]=[x, y]\otimes ab,\ \ x,y\in\lie g, \ \ a,b\in A.$$
\begin{defn}
The Lie algebra $\lie a \otimes A$ is called the \textit{loop algebra} of $\lie a$ and is denoted by $L(\lie a)$.
\end{defn}
We will denote by $\Gamma: A \rightarrow A$ the group action of $\Gamma$ on $A$ given by extending the map $\sigma: t \mapsto \zeta t$ to an algebra homomorphism (recall that $\zeta$ is a primitive $m^{th}$ root of unity).  Then $A$ decomposes as 
\[ A = \bigoplus_{s=0}^{m-1} A_{s}, \]
where $A_s = \left\{ a \in A : \sigma(a) = \zeta^s a \right\}$.  We then have $A_{0}= \bc\left[ t^{\pm m } \right]$ and $A_{s} = t^sA_{0}$.  The linear extension of the map  $\sigma: x \otimes t^k \mapsto \sigma(x) \otimes \sigma(t^k)$ for all $x \in \lie g$, $k \in \bz$ is a Lie algebra automorphism of $\lie g \otimes A$, and the set of fixed points 
\[ (\lie g \otimes A)^\Gamma = \bigoplus_{s =0}^{m-1} \lie g_{s} \otimes A_{-s}\] 
is a Lie subalgebra of $\lie g\otimes A$.
\begin{defn}  The Lie algebra $(\lie g \otimes A)^\Gamma$ is called the \textit{twisted loop algebra} of $\lie g$ with respect to $\Gamma$; we will denote this algebra by $\Lgg$.
\end{defn}
These loop algebras occur as a main ingredient in a realization of the affine Kac--Moody algebras and also of the extended affine Lie algebras; see for example \cite[Chapter 18]{Ca05} or \cite{K90} for details.  For any subalgebra $\lie a$ of $\lie g$ which is invariant under the action of $\Gamma$, we set $L^\Gamma(\lie a) = (\lie a \otimes A)^\Gamma$.

As $\Gamma$ is generated by a diagram automorphism of $\lie g$, the subalgebras $\lie n^\pm, \lie h$ of $\lie g$ are each preserved by $\Gamma$ and hence $L^\Gamma(\lie g)$ inherits the triangular decomposition of $\lie g$:
\[ L^\Gamma(\lie g) = L^\Gamma(\lie n^-) \oplus L^\Gamma(\lie h) \oplus L^\Gamma(\lie n^+).\]

We briefly mention a more geometric realization of these loop algebras.  The ring $A$ is the coordinate ring of the affine variety $\bc^*$.  The Lie algebra $\lie g$ can be viewed as an affine variety, and if we denote by $M(\bc^*, \lie g)$ the Lie algebra of regular maps from $\bc^*$ to $\lie g$ (where the bracket is defined pointwise), the group action of $\Gamma$ on $\lie g$ and on $A$ (hence on $\bc^*$) extends to an action $\Gamma: M(\bc^*, \lie g) \rightarrow M(\bc^*, \lie g)$.  Then it is easy to see that $M(\bc^*, \lie g)^\Gamma \cong L^\Gamma(\lie g)$; we call such a realization of $L^\Gamma(\lie g)$ an \textit{equivariant map algebra} (see \cite{NSS} for more details).  

We identify  $\lie a$ with the  Lie subalgebra  $\lie a\otimes \bc$ of $\lie a \otimes A$.
   Similarly, if $\lie b$ is a Lie subalgebra of $\lie a$, then $\lie b\otimes A$ is
    naturally isomorphic to a subalgebra of $\lie a\otimes A$. Finally we denote by $\bu(\lie g\otimes A_+)$ the subspace of $\bu(\lie g\otimes A)$ spanned by monomials in the elements $x\otimes a$ where $x\in \lie g$, $a\in A_+$.

If $J_{0}$ is any ideal in $A_{0}$, then $\displaystyle{\bigoplus_{s=0}^{m-1} \lie g_s \otimes t^{-s}J_0}$ is clearly an ideal of $L^\Gamma(\lie g)$; conversely, the following can be deduced from \cite{CFS08} or \cite{Lau10}.  

\begin{lem}\label{ideal} 
Let $J$ be an ideal of $L^\Gamma(\lie g)$.  Then there exists an ideal $J_0 \subseteq A_0$ such that $\displaystyle{\mathit{J} = \bigoplus_{s \in \bz} \lie g_{s} \otimes t^{-s}J_{0}}$. 
\end{lem}

\subsection{}
A very important tool for understanding and analyzing modules for loop algebras has been the use of results for  $L(\lie{sl}_2)$. In the twisted loop setting, we will once again use results for the smallest available twisted loop algebra. Namely, let $\lie g = \lie{sl}_3$ and $\Gamma$ be induced by the non-trivial Dynkin diagram automorphism of $\lie g$. Then in our notation, the fixed--point algebra is denoted by $\Lat$. In this case $\lie g_0 \cong \lie{sl}_2$ and $\lie g_1 \cong V(4\omega)$, the five--dimensional irreducible $\lie{sl}_2$-module. \\
In contrast with the loop case, the twisted loop algebras are in some sense built from copies of $L(\lie{sl}_2)$ and $\Lat$. The following lemma, proved in \cite{FK11}, makes this idea precise.

\begin{lem}\label{smallalgebras}\mbox{}
\begin{enumerit}
\item[(i)] If $\lie g$ is of type $A_{2n}$, then we have canonical isomorphisms
 $$L(\lie{sl}_2)\cong\mbox{sp}\left\{x^{\pm}_{\alpha}(j)\otimes t^{ms-j},\  h_{\alpha}(j)\otimes t^{ms-j}\ |\  s\in\mathbb{Z}\ ,\ 0\leq j \leq m-1\right\},$$ if $\alpha$ is a long root, and
 $$\Lat \cong \mbox{sp}\left\{ x^{\pm}_{\alpha}(j)\otimes t^{ms-j},\ x^{\pm}_{2\alpha}(j+1)\otimes t^{ms-j},\  h_{\alpha}(j) \otimes t^{ms-j}\ |\ s\in \mathbb{Z} ,\ 0\leq j \leq m-1\right\},$$
 if $\alpha$ is a short root.
\item[(ii)] If $\lie g$ is not of type $A_{2n}$, then we have canonical isomorphisms
 $$L(\lie{sl}_2)\cong\mbox{sp}\left\{ x^{\pm}_{\alpha}(0) \otimes t^{ms},\ h_{\alpha}(0)\otimes t^{ms}\ |\  s\in\mathbb{Z}\right\},$$
 if $\alpha$ is a long root, and
$$L(\lie{sl}_2)\cong\mbox{sp}\left\{ x^{\pm}_{\alpha}(j)\otimes t^{ms-j},\ h_{\alpha}(j)\otimes t^{ms-j}\ |\  s\in\mathbb{Z}\ ,\ 0\leq j \leq m-1\right\},$$ if $\alpha$ is a short root.
\end{enumerit}

\end{lem}\hfill\qed

\subsection{}
Let $V$ be a $\lie g$--module. We say that $V$ is locally finite--dimensional if
any element  of $ V$ lies in a finite--dimensional $\lie g$--submodule of $V$.  A locally finite--dimensional $\lie g$--module $V$ is isomorphic to a direct
sum of simple finite--dimensional $\lie g$--modules and hence we can write $$V=\bigoplus_{\lambda\in\lie h^*} V_\lambda,$$ where
$V_\lambda=\{v\in V\ :\ h.v=\lambda(h)v,\ \ h\in\lie h\}$. We set
$$\wt(V)=\{\lambda\in\lie h^*:V_\lambda\ne 0\}.$$
 For $\lambda\in P^+$, let $V(\lambda)$ be the finite-dimensional simple
$\lie g$--module of highest weight $\lambda$. It is well known (see, for example, \cite{Hu72}) that $V(\lambda)$ is generated by an element $v_\lambda\in V(\lambda)$ satisfying the defining relations:
$$
\lie n^+.v_\lambda=0,\quad h.v_\lambda=\lambda(h)v_\lambda,\quad (x^-_{i})^{\lambda(h_i)+1}.v_\lambda =0,
$$
for all~$h\in\lie h$, $i\in I$. Moreover, $$\wt(V(\lambda))\subset \lambda -Q^+,\ \ \dim V(\lambda)<\infty,$$ 
and any simple locally finite--dimensional $\lie g$--module is isomorphic to $V(\lambda)$ for some $\lambda\in P^+$.
\begin{rem} We will often use these results in the case for $\lie g_0$, whose simple finite--dimensional modules are then parametrized by elements in $P_0^+$.
\end{rem}

\subsection{}\label{simplemod}
In this section, we will recall results on the classification of simple finite--dimensional modules for $\Lg$. We begin with the definition of an evaluation module. Given $\lambda \in P^+$ and $a \in \bc^*$, the the $\lie g$-module $V(\lambda)$ has an $\Lg$--module structure given by 
$$(x \otimes t^n).v = a^nx.v \ \ \text{ for all }\  x \in \lie g\ , \ v \in V(\lambda).$$
We denote this module by $V_a(\lambda)$. Clearly, since $V(\lambda)$ is a simple $\lie g$-module, $V_a(\lambda)$ is a simple $\Lg$-module. This result has a generalization for tensor products of simple modules.  To state this result, we will first introduce some useful terminology, due to \cite{NSS}:\\
Let $\Xi$ be the monoid of finitely supported functions from $\bc^*$ to $P^+$. Thus, for $\xi \in \Xi$,  
$$\supp(\xi) \coloneqq \{ a \in \bc^* \ | \ \xi(a) \neq 0\} \subset \bc^*$$
is a finite set. We define the \textit{weight} of $\xi\in\Xi$ by the formula $\wt(\xi) \coloneqq \sum_{a \in \supp(\xi)} \xi(a) \in P^+$. Consequently we have
$$\Xi = \bigcup_{\lambda \in P^+} \Xi_\lambda,$$
where $\Xi_\lambda=\{ \xi \in \Xi \ | \ \wt(\xi) = \lambda \}$.
We associate to each $\xi \in \Xi$ an $\Lg$-module
$$V_\xi \coloneqq \bigotimes_{a \in \supp(\xi)} V_a(\xi(a)).$$ The following characterization of simple finite--dimensional $L(\lie g)$-modules was proved in \cite{CP01}.
\begin{thm}
$V_\xi$ is a simple finite--dimensional $\Lg$--module. Moreover, if $V$ is a simple finite--dimensional $\Lg$-module, then there exists $\xi \in\Xi$, such that $V \cong V_\xi$. 
\end{thm}

\subsection{}\label{admissible}
Before recalling the results on simple finite--dimensional modules for $\Lgg$, we will introduce the necessary notion of admissible sets.
\begin{defn} A finite subset $X \subset \bc^*$ is called admissible, if for all $a \neq b \in X$ we have
$$\Gamma.a \cap \Gamma.b = \emptyset.$$
We say a finitely supported function $\xi \in \Xi$ is \textit{admissible} if its support $\supp(\xi)$ is an admissible set. Furthermore, for every finite subset $X \subset  \bc^*$, we denote by $X_{\adm}$ a maximal admissible subset (clearly this set is not unique, but for our purposes the uniqueness will not be necessary).
\end{defn}

Now, observe that any Dynkin diagram automorphism $\sigma$ induces an automorphism of $P^+$ given by the formula $\sigma(\omega_i) = \omega_{\sigma(i)}$. As stated in Section ~\ref{subsectiongeometry}, the group of automorphisms $\Gamma$ acts also on $\bc^*$ by multiplication with $\zeta$, a primitive $m^{th}$ root of unity. \\
We say $\xi\in\Xi$ is \textit{equivariant} with respect to $\Gamma$, if
$$\xi(\sigma(a)) = \sigma(\xi(a)) \ \ \text{for all} \ a \in \bc^* \text{ and } \sigma \in \Gamma.$$
We denote by $\Xi^{\Gamma}$ the set of equivariant functions in $\Xi$.  The following was proved in \cite{NSS}.
\begin{thm}\label{simplelgg}
$\Xi^{\Gamma}$ parametrizes the simple finite--dimensional $\Lgg$-modules.
\end{thm}

For the reader's convenience, we recall here the assigment of a simple module to an equivariant function. In order to do so, we introduce the symmetrizer map $\Sigma: \Xi \longrightarrow \Xi^{\Gamma}$, given by
$$\xi \mapsto \sum_{\sigma \in \Gamma} \sigma \circ \xi \circ \sigma^{-1}.$$
Clearly, this function is well-defined, since the right hand side is by construction equivariant.\\ 
Given $\chi \in \Xi^{\Gamma}$, a function $\xi \in \Xi$ is called $\chi$-\textit{admissible} if $\Sigma(\xi) = \chi$ and $\supp(\xi)$ is an admissible set. Before continuing, we observe that for each $\chi \in \Xi^\Gamma$, there exists at least one $\chi$-admissible function, constructed as follows.

Let  $\chi \in \Xi^{\Gamma}$ and choose a maximal admissible subset $X_{\adm} \subset \supp (\chi) $. For $a\in\bc^*$, define $\xi\in\Xi$ by
$$ \xi(a)\coloneqq \begin{cases} 	\chi(a) & \text{ if } a \in X_{\adm}\\ 	0 & a \notin X_{\adm}\\   \end{cases}.$$
Then $\xi$ is finitely supported, and $\supp (\xi)$ is admissible by construction, with 
$ \Sigma(\chi) = \xi.$

The following was shown in (\cite{Lau10},\cite{NSS}):
\begin{lem}
Suppose $\xi \in \Xi$ is admissible. Then the $\Lg$-module $V_{\xi}$ is simple as an $\Lgg$-module.  Moreover, every simple $\Lgg$-module is obtained in this way.
\end{lem}

The parametrization of Theorem~\ref{simplelgg} is completed by observing that for two admissible functions $\xi_1, \xi_2 \in \Xi$ with $\Sigma(\xi_1) = \Sigma(\xi_2)$, we have
$$V_{\xi_1} \cong V_{\xi_2} \text{ as } \Lgg \text{-modules.}$$

We shall also define the \textit{weight} of an equivariant function $\chi \in \Xi^{\Gamma}$. This was done before for elements from $\Xi$, but it is important to note that although $\Xi^{\Gamma} \subset \Xi$, the weight of an element in $\Xi^{\Gamma}$ considered as an equivariant function is different from its weight considered as an element in $\Xi$. To define the weight of $\chi$, let $\xi \in  \Xi$ be $\chi$-admissible and set
$$\wt_0 (\chi) = \overline{\wt(\xi)} \in  P_0^+.$$
We observe here that the weight is independent of the choice of $\xi$.


\section{The category  \texorpdfstring{$\mathcal{I}^\Gamma $}{of integrable modules}}\label{section3}
In this section we will (by analogy with \cite{CFK10}) define the category of locally finite modules and the global Weyl modules. We keep the exposition as short as possible without sacrificing necessary detail.

\subsection{}
 Let $\mathcal{I}^\Gamma$   be the category whose
objects are modules for $L^\Gamma(\lie g)$ which are locally finite--dimensional $\liefixed g$--modules and whose morphisms are $$\Hom_{\mathcal{I}^\Gamma}(V,V')=\Hom_{L^\Gamma(\lie g)}(V,V'),\ \ V,V'\in\mathcal{I}^\Gamma.$$ Clearly  $\mathcal{I}^\Gamma$ is an abelian category and is closed under tensor products. We shall use the following elementary result often without mention in the rest of the paper.
 \begin{lem}\label{superelem} Let $V\in\Ob \mathcal{I}^{\Gamma}$.  \begin{enumerit} \item If $V_\lambda\ne 0$ and  $\wt V\subset\lambda-Q_0^+$, then $\lambda\in P_0^+$ and $$L^\Gamma(\lie n^+).V_\lambda=0,\ \ (x_i^-)^{\lambda(h_i)+1}.V_\lambda=0,\ \ i\in I_0.$$  If in addition, $V=\bu(L^\Gamma(\lie g)).V_\lambda$ and $\dim V_\lambda=1$, then $V$ has a unique irreducible  quotient.
 \item If $V=\bu(L^\Gamma(\lie g)).V_\lambda$ and $L^\Gamma(\lie n^+).V_\lambda=0$, then $\wt(V)\subset\lambda-Q_0^+$.
 \item If $V\in\mathcal{I}^\Gamma$ is irreducible and finite--dimensional, then there exists $\lambda\in\wt V$ such that 
 $$ \dim V_\lambda=1,\ \ \wt(V)\subset\lambda- Q_0^+.$$
\hfill\qedsymbol
\end{enumerit}
\end{lem}

\subsection{}
We recall here the definition of the global Weyl module for $\Lg$ (due to \cite{CP01}); it will play a crucial role in all that follows. 
\begin{defn}
Let $\lambda \in P^+$.  The global Weyl module $W(\lambda)$ is generated by a non-zero vector $w_\lambda$, subject to the defining relations:
$$ L(\lie n^+).w_\lambda = 0, \quad (H \otimes 1).w_\lambda = \lambda(H)w_\lambda, \quad (X_i^-)^{\lambda(H_i)+1}.w_\lambda = 0, \qquad\  i \in I,\ \  H \in \lie h.$$
\end{defn}
The study of these modules initiated a series of papers (\cite{CP01}, \cite{CL06}, \cite{FoL07}, \cite{Na11}, \cite{BN04}), and we give here a brief summary of the results contained therein. $W(\lambda)$ is an integrable projective module in a certain category (see Section~\ref{subsectioncategory}). Furthermore, $W(\lambda)$ is a free module of finite rank over the algebra 
$$\ba_\lambda\coloneqq\bu(L(\lie h))/\ann_{\bu(L(\lie h))} w_\lambda,$$
which is isomorphic to a Laurent polynomial ring in finitely many variables.

\subsection{}\label{subsectioncategory}
Given an integrable left $\lie g_0$-module $V$, it is a standard fact of relative homological algebra that 
$$P^\Gamma(V) : = \bu(L^\Gamma(\lie g)) \otimes_{\lie g_0} V$$
is a projective $\Lgg$-module.  Moreover $P^\Gamma(V)$ lies in $\mathcal{I}^{\Gamma}$. Furthermore, if $\lambda \in P_0^+$, then $P^\Gamma(V(\lambda))$ is generated as an $L^\Gamma(\lie g)$-module by a non-zero element $v$ with relations
$$\lie n_0^+.v = 0\quad h.v = \lambda(h)v,\quad (x_i^-)^{\lambda(h_i)+1}.v = 0,\qquad \ i\in I_0,\ \ h\in \lie h_0.$$
For $\nu\in P_0^+$ and $V\in\Ob\cal I^\Gamma$, let $V^\nu\in\Ob\cal I^\Gamma$ be defined by: \begin{equation}\label{defglobweylext}V^\nu\coloneqq V/\sum_{\mu\nleq\nu}\bu(L^\Gamma(\lie g))V_\mu. \end{equation} Equivalently, this is the unique  maximal $L^\Gamma(\lie g)$ quotient $W$ of
$V$ satisfying $\wt(W)\subset \nu-Q_0^+.$
A morphism $\pi:V\to V'$ of objects in $\cal I^{\Gamma}$ clearly induces a morphism $\pi^\nu: V^\nu\to (V')^\nu$. Let $\cal I^\Gamma_\nu$ be the full subcategory of objects $V\in\cal I^{\Gamma}$ such that $V=V^\nu$. If $V\in\Ob\cal I^\Gamma_\nu$, then its weights are bounded above by $\nu$ and, since it is integrable, it decomposes into a direct sum simple finite--dimensional $\lie g_0$-modules. Hence, 
\begin{equation}\label{finiteset} V\in\Ob\cal I^\Gamma_\nu\implies \#\wt V<\infty.\end{equation}
\begin{rem}\label{int-weyl} If $\lambda, \nu \in P_0^+$, then $P^\Gamma(V(\lambda))^{\nu} \in \Ob\cal I^\Gamma_\nu$.
\end{rem}
We are now able to define the main object of study for this paper:
\begin{defn}\label{weyldef} The \textit{global Weyl module} of weight $\lambda \in P_0^+$ for $L^\Gamma(\lie g)$ is
$$W^\Gamma(\lambda)\coloneqq P^\Gamma(V(\lambda))^{\lambda}.$$
\end{defn}
The following proposition is proved analogously to \cite[Proposition 4]{CFK10}.
\begin{prop}\label{definingrelations}
$W^{\Gamma}(\lambda)$ is generated by a nonzero element $w_{\lambda}$ with relations
\begin{equation}\label{redweyldef}
L^{\Gamma}(\lie n^+).w_\lambda=0,\quad h.w_\lambda=\lambda(h)w_\lambda,\quad (x^-_{i})^{\lambda(h_i(0))+1}.w_\lambda =0,\quad\ i\in I_0,\ \ h\in\lie h_0.
\end{equation}
\end{prop}

\subsection{}
For $\mu \in P^+$, $W(\mu)$ may be viewed as a module for $\Lgg$ by restriction, in which case the highest weight vector will be of weight $\bar{\mu}$. It follows that there is a natural map $W^{\Gamma}(\bar{\mu}) \longrightarrow W(\mu)$. The immediate questions, whether this map is injective or surjective, must be answered in the negative in general. Nevertheless, we shall find an $\Lg$--module into which a global Weyl module for $\Lgg$ embeds: let $\lambda \in P_0^+$ and consider 
$$W\coloneqq\bigoplus_{\mu \in P^+ \ : \ \bar{\mu} = \lambda} W(\mu).$$
This is clearly a module for $\Lg$ and hence, by restriction, for $\Lgg$. The main result of our paper is that $W^\Gamma(\lambda)$ appears as a submodule of $W$.
\begin{thm}\label{maintheorem}
 Let $\lambda \in P_0^+$.  There is an injective homomorphism of $\Lgg$-modules
$$W^{\Gamma}(\lambda) \hookrightarrow \bigoplus_{\overline{\mu} = \lambda}  W(\mu),$$
induced by the assignment
$$w_{\lambda} \mapsto w \coloneqq \sum_{\overline{\mu} = \lambda} w_\mu.$$
\end{thm}

By Proposition \ref{definingrelations}, it is clear that this assignment gives a homomorphism of $\Lgg$--modules. The remainder of the paper is devoted to the proof that the map is, in fact, injective.


\section{The Weyl functor and its properties}\label{section4}
For $\lambda \in P_0^+$, we denote by $\ann_{\bu(L^\Gamma(\lie h))}w_\lambda$ the annihilator of $w_\lambda$ in $\bu(L^\Gamma(\lie h))$.  This is an ideal in $\bu(L^\Gamma(\lie h))$, and we define $\ba^\Gamma_\lambda$ as the quotient of $\bu(L^\Gamma(\lie h))$ by this ideal:
\[ \ba^\Gamma_\lambda \coloneqq \bu(\Lhg) / \ann_{\bu(\Lhg)}w_\lambda. \]
Clearly, $\Bagg$ is a commutative associative algebra and we will see in Theorem~\ref{fingen} that it is finitely generated.

\subsection{}
We define a right module action of $\ba^\Gamma_\lambda$ on $W^\Gamma(\lambda)$ as follows: for $ a \in \ba^\Gamma_\lambda$ and $u \in \bu(L^\Gamma(\lie g))$, 
\[ uw_\lambda.a \coloneqq ua.w_\lambda.\] 
The verification that this map is well--defined is straightforward; see \cite{CFK10} for details.  For all $\mu \in P_0^+$, the subspaces $W^\Gamma(\lambda)_\mu$ are $\Lhg$--submodules for both the left and right actions and 
\begin{align*} \ann_{\bu(\Lhg)}w_\lambda &= \left\{ u \in \bu(\Lhg) : w_\lambda.u = 0 = u.w_\lambda \right\} \\
&= \left\{ u \in \bu(\Lhg) : W^{\Gamma}(\lambda). u = 0 \right\}. 
\end{align*}
Therefore $W^\Gamma(\lambda)$ is an $\left( L^\Gamma(\lie g), \ba^\Gamma_\lambda\right)$--bimodule and each subspace $W^\Gamma(\lambda)_\mu$ is a right $\ba^\Gamma_\lambda$--module.  Moreover, $W^\Gamma(\lambda)_\lambda$ is an $\ba^\Gamma_\lambda$--module and 
\[ W^\Gamma(\lambda)_\lambda \cong_{\ba^\Gamma_\lambda} \ba^\Gamma_\lambda.\]

Let $\mode \ba^\Gamma_\lambda$ be the category of left $\ba^\Gamma_\lambda$--modules. 
\begin{defn} 
 Let $\bw^\Gamma_\lambda: \mode \ba^\Gamma_\lambda \rightarrow \cal I^\Gamma_\lambda$ be the right exact functor given by 
\[ \bw^\Gamma_\lambda M = W^\Gamma(\lambda) \otimes_{\ba^\Gamma_\lambda} M, \; \; \; \; \; \bw^\Gamma_\lambda f = 1 \otimes f, \]
where $M \in \Ob \mode \ba^\Gamma_\lambda$ and $f \in \Hom_{\ba^\Gamma_\lambda}(M, M')$ for some $M' \in \Ob\mode \ba^\Gamma_\lambda$.  We call this functor the \textit{(twisted) Weyl functor}.
\end{defn}
The $\liefixed g$--action on $\bw^\Gamma_\lambda M$ is locally finite (since $W^\Gamma(\lambda) \in \Ob\cal I^\Gamma_\lambda$), so that $\bw^\Gamma_\lambda M \in \Ob \cal I^\Gamma_\lambda$, and 
\[ \bw^\Gamma_\lambda \ba^\Gamma_\lambda \cong_{L^\Gamma(\lie g)} W^\Gamma(\lambda), \; \; \; \; \; (\bw^\Gamma_\lambda M)_\mu \cong_{\ba^\Gamma_\lambda} W^\Gamma(\lambda)_\mu \otimes_{\ba^\Gamma_\lambda} M,\]
for $\mu \in P_0$, $M \in \Ob\mode \ba^\Gamma_\lambda$.

\subsection{}
In this section we adapt results from \cite{CFK10} and state them without proofs, since these proofs carry over almost verbatim from the case of untwisted loop algebras. The first ingredient we shall require is the restriction functor, which will be right adjoint to the Weyl functor. 
For this we need the following lemma, whose proof can be found in \cite[Lemma 4]{CFK10}.
\begin{lem}\label{hwmodules}
For all $\lambda \in P_0^+$ and $V \in  \Ob \cal I^{\Gamma}_\lambda$ we have $\ann_{L^{\Gamma}(\lie h)} (w_\lambda) . V_\lambda = 0$.  
\end{lem}

As a consequence, we see that the left action of $\bu(\Lhg)$ on $V_\lambda$ induces a left action of $\ba^\Gamma_\lambda$ on $V_\lambda$, and we denote the resulting $\ba^\Gamma_\lambda$--module by $\br^\lambda_\Gamma V$.  Given
$\pi \in \Hom_{\cal I^{\Gamma}_\lambda}(V,V')$ the restriction $\pi_\lambda:V_\lambda\to V'_\lambda$ is a morphism of $\ba^\Gamma_\lambda$--modules and
$$V\mapsto\br^\lambda_\Gamma V,\ \ \pi\mapsto \br^\lambda_\Gamma \pi=\pi_\lambda $$ 
defines a   functor 
$$\br^\lambda_\Gamma:\cal I^\Gamma_\lambda\longrightarrow \mode\ba^\Gamma_\lambda$$ 
which is exact since restriction of $\pi$ to a weight space is exact.  If $M\in\Ob\mode\ba^\Gamma_\lambda$, we have an isomorphism of left $\ba^\Gamma_\lambda$--modules, 
$$\br^\lambda_\Gamma\bw^\Gamma_\lambda M=(\bw^\Gamma_\lambda M)_\lambda=W^{\Gamma}(\lambda)_\lambda\otimes_{\ba^\Gamma_\lambda}M\cong w_\lambda\ba^\Gamma_\lambda\otimes_{\ba^\Gamma_\lambda} M\cong M,$$
and hence an isomorphism of functors $\id_{\ba^\Gamma_\lambda}\cong \br^\lambda_\Gamma\bw^\Gamma_\lambda$.\\

We may apply the restriction functor to an object of $\cal I^\Gamma_\lambda$, then apply the Weyl functor and obtain once again an object of $\cal I^\Gamma_\lambda$. The following proposition shows the relationship between these two $\Lgg$--modules.
\begin{prop}\label{rightadjoint}
Let $\lambda\in P_0^+$ and  $V\in \Ob \cal I^{\Gamma}_\lambda$. There exists a canonical map of $L^\Gamma(\lie g)$--modules $\eta_V: \bw^\Gamma_\lambda\br^\lambda_\Gamma V\to V$ such that $\eta:\bw^\Gamma_\lambda\br^\lambda_\Gamma\Rightarrow\id_{\cal I^{\Gamma}_\lambda}$ is a natural transformation of functors and $\br^\lambda_\Gamma$ is a right adjoint to $\bw^\Gamma_\lambda$.
\end{prop}
As an immediate consequence we obtain with standard homological methods:
\begin{cor}
 The functor $\bw^\Gamma_\lambda$ maps projective objects to projective objects.
\end{cor}
We now have all the ingredients necessary to state the main result (Theorem 1) of \cite{CFK10} in the case of twisted loop algebras. The proof carries over almost identically, hence will be omitted.

\begin{thm}
Let $\lambda \in P_0^+$ and $V\in\Ob\cal I_\lambda^{\Gamma}$. Then $V\cong\bw^\Gamma_\lambda\br^\lambda_\Gamma V$ iff for all $U\in\Ob \cal I^{\Gamma}_\lambda$ with $U_\lambda=0$, we have \begin{equation}\label{altweyldef} \Hom_{\cal I^{\Gamma}_\lambda}(V,U)=0,\ \ \Ext^1_{\cal I^{\Gamma}_\lambda}(V, U)=0.\end{equation}
\end{thm}

\subsection{}
Another consequence of Proposition~\ref{rightadjoint} is the one-to-one correspondence between maximal ideals of $\Bagg$ and simple modules of $\Lgg$ of highest weight $\lambda$. 
\begin{lem}\label{baggmax}
For $\lambda \in P_0^+$, there exists a natural correspondence between maximal ideals of $\Bagg$ and $\Xi_\lambda^{\Gamma}$.
\end{lem}
\begin{pf}
Let $\bm \in \Max \Bagg$. Then $\bw_\lambda^{\Gamma} (\Bagg/\bm)$ has a unique simple quotient of highest weight $\lambda$, denoted by $V_{\xi_\bm}$ (as in Section~\ref{simplemod}). On the other hand, let $\xi \in \Xi_\lambda^{\Gamma}$ and $V_{\xi}$ be the corresponding simple $\Lgg$-module (Theorem~\ref{simplelgg}).  Then $\br^\lambda_\Gamma V_\xi$ is a simple $\Bagg$-module, and so there exists $\bm_\xi \in \Max \Bagg$ such that 
$$\br^\lambda_\Gamma V_\xi \cong \Bagg/\bm_{\xi}.$$
Since $\br^\lambda_\Gamma$ is right adjoint to $\bw_\lambda^{\Gamma}$, we have $\xi_{\bm_\xi} = \xi$ and $\bm = \bm_{\xi_\bm}$.
\end{pf}

\subsection{}
In this section, we will prove that the global Weyl module is finitely generated as a right $\ba_\lambda^\Gamma$--module. This result is analogous to \cite[Theorem 2]{CFK10}, but requires a new proof when $\Gamma \neq \Id$. \\
To clarify the importance of this result, first recall that the global Weyl module is infinite--dimensional, and even decomposes into infinitely many simple $\lie g_0$-modules. 
By applying the Weyl functor on one--dimensional $\Bagg$-modules we obtain the so called \textit{local Weyl modules} (see Section~\ref{local-comp}). Once we show that the global Weyl module is finitely generated, we can deduce that the local Weyl modules are finite--dimensional (see \cite[Proposition 4.2]{CFS08}). 

\begin{thm} \label{fingen}
$W^{\Gamma}(\lambda)$ is a finitely generated right $\ba^\Gamma_\lambda$--module.
\end{thm}

Let $u$ be an indeterminate and for $\alpha \in R_0^+$, we define for all $\ell \geq 1$ the following power series in $u$ with coefficients in $\bu(\Lhg)$: 
\[ \bop_{\alpha,\ell}(u) =  \exp\left(  - \sum_{k=1}^\infty  \frac{h_\alpha \otimes t^{\ell k}}{k}u^{k} \right). \] 
Let $\bop_{\alpha,\ell}(u) = \sum_{j=0}^\infty p_{\alpha,\ell}^j u^j$.  Note that $p_{\alpha,\ell}^0 = 1$, and that $p_{\alpha,\ell}^j$ is contained in the subalgebra generated by $\left\{ h_\alpha \otimes t^{\ell k} : 0 \leq k \leq j \right\}$.
For the proof, we need the following lemma, which is an immediate consequence of \cite[Lemma 2.3]{CFS08} via the substitution $t \mapsto t^\ell$.  

\begin{lem}\label{Garland} Let $\alpha \in R_0^+$ and $r\in\bz_+$.  Then if
$$ \ell \in \begin{cases} \mathbb{Z}_{> 0} \ , \ \text{ if } \alpha \in (R_0)_s \text{ and } \lie g \text{ is not of type } A_{2n}\\
m\mathbb{Z}_{> 0} \ , \ \text{ if } \alpha \in (R_0)_l \text{ and } \lie g \text{ is not of type } A_{2n}\\
2\mathbb{Z}_{> 0} \ , \ \text{ if } \alpha \in (R_0)_s \text{ and } \lie g \text{ is of type } A_{2n}\\
\mathbb{Z}_{> 0} \ , \ \text{ if } \alpha \in (R_0)_l \text{ and } \lie g \text{ is of type } A_{2n}
\end{cases}
$$
we have
$$ (x_\alpha^+\otimes t^\ell)^r \left( x_\alpha^- \otimes 1 \right)^{r+1}
+(-1)^{r+1} \left( \sum_{j=0}^{r} \left( x_\alpha^-\otimes t^{\ell j}\right) p_{\alpha,\ell}^{r-j} \right) \in \bu(L^\Gamma(\lie g))\bu(L^\Gamma(\lie n^+))_+.$$
\end{lem}\hfill\qedsymbol

Using this containment, we now prove the theorem.
\begin{pf}
Since $W^\Gamma(\lambda)$ is an object of $\mathcal{I}_\lambda^{\Gamma}$ (see Remark~\ref{int-weyl}), we know that $\#\wt(W^\Gamma(\lambda)) < \infty$. It follows that for any monomial $u \in \bu(\Lgg)$ with $\wt_{\lie g_0} (u) \neq 0$, and any $v \in W^{\Gamma}(\lambda)$, there exists $N > 0$ such that $u^N.v = 0$. In particular, $W^{\Gamma}(\lambda)$ is locally finite--dimensional for any vector of the form $x \otimes t^k$, $x \in \lie n^- \cap \lie g_{\epsilon}$, $k \equiv -\epsilon \mod m$.\\

Let $\lambda \in P_0^+$ and $\alpha \in R_0^+$. We set $\ell = 1$ if $\lie g$ is of type $A_{2n}$ and $\alpha$ is a long root, or $\lie g$ is not of type $A_{2n}$ and $\alpha$ is a short root. In any other case, we set $\ell = m$. By setting $r = \lambda(h_\alpha)$, we see from the defining relations of $W^\Gamma(\lambda)$ and the lemma above that 

\begin{equation}\label{garlandcor2}
 x_\alpha^-\otimes t^{\ell r} .w_\lambda = (-1)^r\left( \sum_{j=0}^{r-1}x_\alpha^- \otimes t^{\ell j} p_{\alpha,\ell}^{r-j} \right).w_\lambda,\end{equation}

which, after an inductive argument, implies
\[ \left( x_\alpha^- \otimes t^{\ell k} \right).w_\lambda \in {\mbox{sp}} \left\{ \left( x_\alpha^- \otimes t^{\ell s} \right)w_\lambda \ba^\Gamma_\lambda : 0 \leq s < \lambda(h_\alpha) \right\}.\]

Additionally we must consider the elements $x^-_{\nu} \otimes t$ and $x^-_{2 \nu} \otimes t$, when $\lie g$ is of type $A_{2n}$ and $\nu \in (R_0)_s$. We proceed with the latter case, the former being very similar.\\
Set $\beta = 2 \nu$ and let $\lie a$ be the Lie algebra generated by the elements
\[ x^+_\beta \otimes t^{2q + 1}, \qquad x^-_\beta \otimes t^{2q-1}, \qquad \frac{1}{2}h_\nu \otimes t^{2q}, \qquad q \in \bz, \]
then we have a Lie algebra ismorphism to $L(\lie{sl}_2)$ (by Lemma~\ref{smallalgebras}), given by
$$ x^+_\beta \otimes t^{2q + 1} \mapsto x_{\alpha}^+ \otimes t^q \ , \ x^-_\beta \otimes t^{2q-1} \mapsto x_\alpha^- \otimes t^q \ , \ \frac{1}{2}h_\nu \otimes t^{2q} \mapsto h \otimes t^q.$$
Lemma~\ref{Garland} gives us 
\[ (x_\beta^+\otimes t^3)^r \left( x_\beta^- \otimes t^{-1} \right)^{r+1}
+(-1)^{r+1} \left( \sum_{j=0}^{r} \left( x_\beta^-\otimes t^{2j-1}\right) p_\beta^{r-j} \right) \in \bu(L^\Gamma(\lie g))\bu(L^{\Gamma}(\lie n^+)_+),\]
where we define $p_\beta^j \in \ba_\lambda^\Gamma$ by 
\[ \sum_{j=0}^\infty p_\beta^j u^j = \exp \left( - \sum_{k=1}^\infty \frac{\frac{1}{2} h_\nu \otimes t^{2k}}{k}u^k\right).\]

Again, since $W^{\Gamma}(\lambda)$ is integrable, we have 
$$x_{\beta}^{-} \otimes t^{2r-1}.w_{\lambda} = (-1)^r\left( \sum_{j=0}^{r-1} \left( x_\beta^-\otimes t^{2j-1}\right) p_\beta^{r-j} \right).w_{\lambda}$$
for $r \gg 0$, and the second case is proven.

To complete the proof of the theorem, let $\left\{ \beta_1, \cdots, \beta_N \right\}$ be an enumeration of $R_0^- \cup R_1^-$ (resp. $R_0^- \cup R_1^- \cup R_2^-$). Using the PBW theorem, we can see that elements of the form 
\[  \left( \left( \lie g_{\epsilon_1}\right)_{\beta_{i_1}} \otimes t^{r_1}\right) 
\left( \left(\lie g_{\epsilon_2}\right)_{\beta_{i_2}} \otimes t^{r_2}\right)  \cdots \left( \left( \lie g_{\epsilon_\ell}\right)_{\beta_{i_\ell}} \otimes t^{r_\ell}\right).w_\lambda   \]
for $0 \leq \epsilon_j < m$, $\beta_{i_j} \in \left\{ \beta_1, \ldots, \beta_N \right\}$, $1 \leq i_1 \leq i_2 \leq \cdots \leq i_\ell \leq N$, $\ell  \in \bz_+$, and $r_i \equiv \epsilon_i \mod m$, generate $W^\Gamma(\lambda)$ as a right module for $\ba^\Gamma_\lambda$. Using this spanning set and the fact that $\#\wt(W^\Gamma(\lambda)) < \infty$, an inductive argument on the length $\ell$ of a monomial from this spanning set shows that $W^\Gamma(\lambda)$ is finitely generated as an $\ba^\Gamma_\lambda$--module.  
\end{pf}

\section{The algebra  \texorpdfstring{$\Bagg$}{of the highest weight space}}\label{section-bagg}
In this section we will give an explicit description of the algebra $\ba^\Gamma_\lambda$ and deduce that it is finitely generated. 

\subsection{}\label{sym-def}
We continue to follow the model of \cite{CFK10}. That is, we identify $\ba^\Gamma_\lambda$ with a ring of symmetric polynomials. To begin, recall that a basis for $\Lhg$ is given by the set $$\{h_i(\epsilon)\otimes t^{mk-\epsilon}\ :\ i\in I_0,\ 1\le \epsilon\le m-1,\ k\in\bz\}=\{h_i(\bar k)\otimes t^{-k}\ :\ k\in \bz,\ i\in I_0\},$$  
where $\bar k$ denotes the least nonnegative residue of $k$ modulo $m$.

Set $A\left( k \right) = \bc\left[ t^{\pm k} \right]$ and for $\lambda = \sum_{i \in I_0}r_i \omega_i \in P_0^+$, define 
\[ \mathbb A^\Gamma_\lambda = \bigotimes_{i \in I_0} \left( A( \left| \Gamma_i \right|)^{\otimes r_i} \right)^{S_{r_i}}. \]

Now we identify a natural generating set of $\mathbb A^\Gamma_\lambda$. 

For $N\in\bz_+$ and $B$ any associative algebra over $\bc$, we define a homomorphism of algebras $\sym_N: B\to B^{\otimes N}$ by the assignment \[ b\mapsto \sum_{\ell=0}^{N-1}1^{\otimes \ell}\otimes b\otimes 1^{\otimes N-\ell-1}.\] 
Now set  \[\sym_{\lambda}^i(b)=1^{\otimes \lambda(\sum_{j<i} h_j)}\otimes\sym_{\lambda(h_i)}(b)\otimes1^{\otimes\lambda(\sum_{j>i} h_j)}\] 
for $i\in I_0$. Taking $B=A(|\Gamma_i|)$ for $i\in I_0$, we clearly have  \[\{\sym_\lambda^i(t^{k})\ :\  k\in|\Gamma_i|\bz\}\subset \mathbb A^\Gamma_\lambda.\] The following lemma makes clear the role of these elements in generating $\mathbb A^\Gamma_\lambda$.

\begin{lem}\label{fancy-genset} The set $$\{\sym_\lambda^i(t^{k})\ :\ i\in I_0,\ k\in|\Gamma_i|\bz,\  |k|\le\lambda(h_i)\}$$ generates $\mathbb A^\Gamma_\lambda$ as an algebra over $\bc$. \end{lem}
\begin{pf}
Fix $i\in I_0$ and set $N=\lambda(h_i)$. It is well known that the algebra $A((|\Gamma_i|)^{\otimes N})^{S_N}$ is isomorphic to the polynomial algebra $\bc[f_1,f_2\ldots,f_N,f_N^{-1}]$, where the $f_\ell$ are the elementary symmetric functions in the $N$ variables $t_1^{|\Gamma_i|},\ldots,t_N^{|\Gamma_i|}.$  On the other hand, the element $\sym_\lambda^i(t^k)$ corresponds to the power sum $g_k$ of degree $k|\Gamma_i|$ in the factor $A(|\Gamma_i|)^{\otimes N}$, so that we have \[ \bc[g_1,\ldots,g_N]\cong\bc[t_1^{|\Gamma_i|},\ldots,t_N^{|\Gamma_i|}]^{S_N}\cong\bc[f_1,\ldots,f_N].\] In order to prove the lemma, it therefore suffices to check that $f_N^{-1}$ lies in $\bc[g_{-1},\ldots,g_{-N}]$. But this is clear, since \[f_N^{-1}\in\bc[t_1^{-|\Gamma_i|},\ldots,t_N^{-|\Gamma_i|}]^{S_N}\cong\bc[g_{-1},\ldots,g_{-N}].\] 
\end{pf}

As a consequence of Lemma~\ref{fancy-genset}, we see that the assignment $$h_i(\bar k)\otimes t^{-k}\mapsto \sym_{\lambda}^i(t^{-k}),\ \ i\in I_0,\ \ k\in\bz$$ extends to a surjective homomorphism of algebras $\tilde\tau_{\lambda}: \bu(\Lhg)\twoheadrightarrow \mathbb A^\Gamma_\lambda$. We shall show in the rest of this section that $\tilde\tau_\lambda$ descends to an isomorphism $\tau_\lambda: \Bagg\cong\mathbb A^\Gamma_\lambda$.

\subsection{} 

The first step is to provide a natural correspondence between $\Xi^\Gamma_\lambda$ and the maximal spectrum of $\mathbb A^\Gamma_\lambda$. This description of $\Max(\mathbb A^\Gamma_\lambda)$ will be used in the sequel to show that $\mathbb A_\lambda^\Gamma$ and $\ba^\Gamma_\lambda$ are isomorphic as algebras.  To describe the correspondence, we introduce an alternate description of the 
maximal ideals in $\mathbb A^\Gamma_\lambda$ in terms of multisets on the maximal ideals of $A(\left| \Gamma_i \right|)$. 

For any set $S$, let $\mathcal{M}(S)$ be the set of functions $f: S\to \bz_+$ satisfying the condition that $f(s)=0$ for all but finitely many $s\in S$. Such a function is called a \textit{finite multiset} on $S$. $\cal M(S)$ forms a commutative monoid under the usual addition of functions. The size of $f\in \cal M(s)$ is given by the formula\[ |f|=\sum_{s\in S} f(s).\] We also note that any element of $\cal M(S)$ can be written uniquely as a $\bz_+$-linear combination of characteristic functions $\chi_s$ for $s\in S$, defined by $\chi_s(b) = \delta_{s, b}$ for $b\in S$. 
\subsection{}
We shall use this language to describe the maximal ideals of the tensor product $\mathbb A^\Gamma_\lambda$. It is clearly enough to classify the maximal spectrum of rings  $(A(|\Gamma_i|)^{\otimes \ell})^{S_\ell}$ for $\ell\in\bz_+$. Such ideals are given precisely by unordered $\ell$-tuples (with possible repetition) of maximal ideals of $A(|\Gamma_i|)$, $i\in I_0$, i.e. by elements $f\in\mathcal{M}\left( \Max(A( \left| \Gamma_i \right|  ))\right)$ with $|f|=\ell$.

Since the maximal ideals of $A(|\Gamma_i|)$ are principal ideals generated by polynomials $t^{|\Gamma_i|}-a^{|\Gamma_i|}$ for $a\in\bc^*$, we may view elements of $\mathcal{M}\left( \Max(A( \left| \Gamma_i \right|  ))\right)$ as multisets consisting of orbits of $\bc^*$ under the action of $\Gamma_i$.

Abbreviating $\mathcal{M}\left( \Max(A( \left| \Gamma_i \right|  ))\right)$ by $\mathcal{M}_i$, the product $\hat{\cal M}=\prod_{i \in I_0}\mathcal{M}_i$ is also a commutative monoid, and
for $\hat f\in\hat{\cal M}$ we set $\wt(\hat f)=\sum_{i\in I_0}|f_i|\omega_i\in P_0^+.$ Defining 
\[ \hat{\mathcal{M}}_\lambda = \left\{ \hat f \in \hat{\cal{M}} : \wt(\hat f)=\lambda \right\},\quad \lambda\in P_0^+,  \]  we see that $\Max\left(\mathbb A_\lambda^\Gamma\right)$ is in bijective correspondence with $\hat{\mathcal{M}}_\lambda.$

In this language, any  $\xi\in\Xi$ can be written uniquely as a linear combination of fundamental weights $\omega_i$, $i \in I$, with coefficients from $\mathcal{M}(\bc^*)$: 
\[ \Xi^\Gamma = \left\{ \sum_{i \in I} f_i \omega_i : f_i \in \mathcal{M}(\bc^*) \right\},\]
where for $\xi = \sum_{i \in I} f_i \omega_i$ and $c \in \Max(A)$, we have $\xi(c) = \sum_{i \in I} f_i(c) \omega_i \in P^+$.

\subsection{}
Next, we exhibit an isomorphism of monoids between $ \hat{\cal{M}}$ and $\Xi^\Gamma$.  Observe that
 for each $i\in I_0$ there is a surjective morphism of monoids
\[ \mathcal{M}(\bc^*) \stackrel{\pi_i}{\twoheadrightarrow} \mathcal{M}_i, \]
defined by extending the assignment $\chi_a\mapsto\chi_{\bar{a}}$, where $\bar{a}$ is the orbit of $a$ under the action of $\Gamma_i$. If $\Gamma_i$ is trivial, then of course $\pi_i$ is simply the identity map. We describe some of its further properties in the following lemma.

\begin{lem}\label{pi-prop} Let $\xi=\sum_{i\in I} f_i\omega_i$ lie in $\Xi^\Gamma$. Then 
\begin{enumerate}
\item[(1)] For all $\gamma\in \Gamma$ and $i\in I$ we have $\pi_i(f_i) = \pi_i(f_{\gamma(i)})$.
\item[(2)] For each $a\in\bc^*$ and $i\in I$, we have 
\begin{equation}\label{pi-f}\pi_i(f_i)(\bar a)=|\Gamma_i|f_i(a)\in |\Gamma_i|\bz_+.\end{equation} 
\end{enumerate}
\end{lem}
\begin{pf}
For the first part, observe that by the equivariance of $\xi$ we have \begin{equation}\label{equizeta}f_i(\gamma(a))=f_{\gamma(i)}(a)\end{equation} for $i\in I$, $\gamma\in\Gamma$ and $a\in 
\Max(A)$. The result now follows by the definition of $\pi_i$.

For the second assertion, there is nothing to prove unless $\Gamma_i=\Gamma$, that is $\gamma(i)=i$ for all $\gamma\in\Gamma$, in which case the result follows from Equation \eqref{equizeta} and the definition of $\pi_i$.
\end{pf}

Now construct a morphism of monoids $\alpha: \Xi^\Gamma\to \hat{\mathcal{M}}$: given an equivariant function $\xi \in \Xi^\Gamma$, write 
$\xi = \displaystyle\sum_{i \in I} f_i \omega_i$ and define $\alpha(\xi) \in \hat{\cal M}$ by the formula
\begin{equation}\label{bij2} \alpha(\xi) \coloneqq \left(\frac{1}{|\Gamma_i|} \pi_i(f_i) \right)_{i \in I_0}.\end{equation} 
It follows from part (ii) of Lemma~\ref{pi-prop} that $\alpha$ is injective, so it remains to show that it is surjective. For this, fix \[\hat g=(g_i)_{i\in I_0}\in \hat{\cal M}.\] To find a preimage of $\hat g$ under $\alpha$, define $f_i\in\cal M(\bc^*)$ for $i\in I_0$ by \[ f_i(a)=\left\{\begin{array}{ccc} g_i(a)&if&\Gamma_i=1,\\ g_i(\bar{a})&if&\Gamma_i=\Gamma,\end{array}\right.\] so that $f_i$ is clearly constant on the orbits of $\bc^*$ under $\Gamma_i$. It follows that $f_i$ satisfies Equation~\ref{pi-f}, and hence $\pi_i(f_i)=|\Gamma_i|g_i$. Hence, $\alpha$ is an isomorphism.

Finally, we show that $\alpha$ induces a bijection $\Max(\mathbb A_\lambda^\Gamma) \leftrightarrow \Xi_\lambda^\Gamma$. For this, it suffices to show that for $\xi\in\Xi^\Gamma$, we have
\begin{equation}\label{weight-zero} \wt_0(\xi)=\wt(\alpha(\xi)),\end{equation} which follows from the observation that \[\wt_0(\xi) = \sum_{i \in I_0} \frac{1}{\left| \Gamma_i \right|} \left| f_i \right| \omega_i,\quad \text{for}\ \  \xi=\sum_{i\in I} f_i\omega_i.\]

\subsection{}\label{taudescends}
The next step is to show that $\tilde{\tau}_{\lambda}$ descends to a surjective homomorphism of algebras \[ \tau_{\lambda}: \mathbf A_{\lambda}^\Gamma\to \mathbb A_{\lambda}^\Gamma.\] 
 
For $\xi\in\Xi^\Gamma_\lambda$, define $\eval_\xi: \bu(\Lhg)\to\bc$ by extending the assignment \begin{equation}\label{}h_i(\bar k)\otimes t^{-k}\mapsto \sum_{\bar{a}\subset \supp(\xi)} a^{-k} \xi(a)(h_i(\bar k)),\quad i\in I_0,\ \ k\in\bz,\end{equation}
 so that \begin{equation*} u.v=\eval_\xi(u)v,\ \ \ v\in(V_\xi)_\lambda,\ \ u\in\bu(L^\Gamma(h)).\end{equation*} 
 Since $V_\xi$ is a quotient of $W^\Gamma(\lambda)$ for $\xi\in\Xi_\lambda^\Gamma$, it follows immediately that \begin{equation}\label{anncont}\ann_{\bu(\Lhg)}w_\lambda\subseteq \bigcap_{\xi\in\Xi_\lambda^\Gamma} \ker(\eval_\xi).\end{equation}

For $\hat f\in\hat{\cal M}_\lambda$, write $\eval_{\hat f}:\mathbb A_\lambda^\Gamma\to \bc$ for evaluation at the maximal ideal of $\mathbb A_\lambda^\Gamma$ corresponding to $\hat f$. Applying the relevant definitions, we have

\begin{equation}\label{commdiag}
\eval_{\alpha^{-1}(\hat f)}=\eval_{\hat f}\circ \tilde{\tau}_\lambda.
\end{equation}

We can now complete the proof that $\mathbb A_\lambda^\Gamma$ is a quotient of $\mathbf A_\lambda^\Gamma$.
\begin{prop}
\[ \ann_{\bu(\Lhg)} w_{\lambda}^\Gamma\subseteq \ker(\tilde{\tau}_{\lambda}).\]
\end{prop}
\begin{pf}
It follows immediately from Equation \ref{commdiag} that 
\begin{equation}\label{kernelint} \bigcap_{\hat f\in\hat{\cal M}_\lambda}\ker(\eval_{\hat f}\circ \tilde{\tau}_\lambda)=\bigcap_{\xi\in{\Xi_\lambda^\Gamma}}\ker(\eval_\xi).\end{equation}
On the other hand, since the Jacobson radical $J(A(s))=0$ for all $s\in\bz_+$, we see that $J(\mathbb A_\lambda^\Gamma)=0$. In particular, 
\begin{equation}\label{jrad} \ker(\tilde{\tau}_\lambda)=\bigcap_{\hat f\in\hat{\cal M}_\lambda}\ker(\eval_{\hat f}\circ \tilde{\tau}_\lambda).\end{equation}
The proposition now follows from Equations \ref{anncont}, \ref{kernelint} and \ref{jrad}.
\end{pf}
\begin{cor}
The map $\tilde{\tau}_{\lambda}$ descends to a surjective homomorphism of algebras \[ \tau_\lambda: \mathbf A_\lambda^\Gamma\to \mathbb A_\lambda^\Gamma.\]
\end{cor}

\subsection{}
It remains to show that $\tau_{\lambda}$ is injective. For this, we adapt the argument of \cite{CFK10}, by identifying a spanning set of $\Bagg$ which is mapped to a linearly independent subset of $\mathbb A^\Gamma_{\lambda}$. 

\begin{lem}\label{bagg-span}
The images of elements \[\left\{ \prod_{i\in I_0}\prod_{s=1}^{m_i}(h_i(\bar{k}_{i,s})\otimes t^{-k_{i,s}})\ :\ 0\le m_i\le \lambda(h_i),\ k_{i,s}\in\bz\right\}\] span ${\Bagg}$.
\end{lem}
\begin{pf}
It is clearly enough to prove that for each $i\in I_0$ and $k_1,\ldots,k_{\ell}\in\bz$ we have \begin{equation}\label{h-span}\prod_{s=1}^\ell(h_i(\bar{k_s})\otimes t^{-k_s})w_\lambda\in{\text{sp}}\left\{ \prod_{s=1}^r (h_i(\bar{\ell_s})\otimes t^{-\ell_s})w_\lambda\ :\ r\le \lambda(h_i)\right\}.\end{equation}
We shall prove this statement as Corollary~\ref{h-cor} below. Assuming it, the lemma follows.
\end{pf}

In order to establish Equation~\ref{h-span}, we shall first prove the following, more general, proposition. For any element $a$ of an associative algebra and any $n\in \bz$, we denote by $a^{(n)}$ the divided power $a^n/n!$, with the convention that $a^{(n)}=0$ for $n<0$.

\begin{prop}\label{commute}
Let $k, \ell\in\bz_+$ with $k\le \ell$. Given $\epsilon_1,\ldots,\epsilon_k\in\bz$ and $i\in I_0$, we have $$ \prod_{s=1}^k(x_i^+(\bar\epsilon_s)\otimes t^{-\epsilon_s})(x_i^-)^{(\ell)}
=\sum_{r=0}^{2k}(x_i^-)^{(\ell-r)} P_r(\epsilon_1,\ldots,\epsilon_k),$$ where $P_r(\epsilon_1,\ldots,\epsilon_k)$ is an element of $U(\Lgg)$ in the standard PBW order, having length at most $k$ and consisting of homogeneous elements of weight $(k-r)\alpha_i$ as an $\lie{sl}_2(i)$-module.
Moreover, $P_k(\epsilon_1,\ldots,\epsilon_k)$ has, except for terms ending in $L^\Gamma(\lie n^+)$, a unique term of length $k$, which is $\prod_{s=1}^k(h_i(\bar\epsilon_s)\otimes t^{-\epsilon_s}).$
\end{prop}
\begin{pf}
The proof proceeds by induction on $k$. For the base case, a simple induction on $\ell$ shows that 
\begin{align}\label{indbase} (x^+_i(\bar{\epsilon})\otimes
 t^{-\epsilon})(x_i^-)^{(\ell)}&=(x_i^-)^{(\ell)}(x_i^+(\bar{\epsilon})\otimes
 t^{-\epsilon}) \\   &+(x_i^-)^{(\ell-1)}(h_i(\bar{\epsilon})\otimes
 t^{-\epsilon})+(x_i^-)^{(\ell-2)}(-x^-_i(\bar\epsilon)\otimes t^{-\epsilon}). \notag\end{align}
Now assuming the result for $k<\ell$, we prove it for $k+1$. By the induction hypothesis and repeated use of Equation \ref{indbase}, we have 
\begin{align*}\prod_{s=1}^{k+1}(x_i^+(\bar\epsilon_s)\otimes t^{-\epsilon_s})(x_i^-)^{(\ell)}&=
\sum_{r=0}^{2k}(x_i^-)^{(\ell-r)}(x_i^+(\bar\epsilon_{k+1})\otimes t^{-\epsilon_{k+1}})P_r(\epsilon_1,\ldots,\epsilon_k)\\
 &+ \sum_{r=0}^{2k}(x_i^-)^{(\ell-r-1)}(h_i(\bar\epsilon_{k+1})\otimes t^{-\epsilon_{k+1}})P_r(\epsilon_1,\ldots,\epsilon_k).\\
&+\sum_{r=0}^{2k}(x_i^-)^{(\ell-r-2)}(x_i^-(\bar\epsilon_{k+1})\otimes t^{-\epsilon_{k+1}})P_r(\epsilon_1,\ldots,\epsilon_k).\end{align*}
Reindexing, this is 
\[ \sum_{r=0}^{2(k+1)}(x_i^-)^{(\ell-r)}P_r(\epsilon_1,\ldots,\epsilon_{k+1}),\]
where 
\begin{align*} P_r(\epsilon_1,\ldots,\epsilon_{k+1})&=(x_i^+(\bar\epsilon_{k+1})\otimes t^{-\epsilon_{k+1}})P_r(\epsilon_1,\ldots,\epsilon_k)\\
&+ (h_i(\bar\epsilon_{k+1})\otimes t^{-\epsilon_{k+1}})P_{r-1}(\epsilon_1,\ldots,\epsilon_{k+1})
+ (x_i^-(\bar\epsilon_{k+1})\otimes t^{-\epsilon_{k+1}})P_{r-2}(\epsilon_1,\ldots,\epsilon_k).\end{align*} 
(Here $P_j(\epsilon_1,\ldots,\epsilon_k)=0$ if $j<0$ or $j>2k$.)
This element, once it is commuted into PBW order, clearly has the correct weight and maximum length. 

It only remains to analyze $P_{k+1}(\epsilon_1,\ldots,\epsilon_{k+1}).$ Any monomial from the third term of this element ends in a summand of $P_{k-1}(\epsilon_1,\ldots,\epsilon_k)$, which has weight $2$ and hence, being already in PBW order, must end in some term from $L^\Gamma(\lie n^+)$. By the induction hypothesis, the second term contains the desired product \[\prod_{s=1}^{k+1} (h_i(\bar\epsilon_{s})\otimes t^{-\epsilon_{s}})\] as its unique term of length $k+1$ (modulo $L^\Gamma(\lie n^+)$).

To deal with the term \[ (x_i^+(\bar{\epsilon}_{k+1})\otimes t^{-\epsilon_{k+1}})P_r(\epsilon_1,\ldots,\epsilon_k),\] we observe that by weight considerations, any monomial not ending in $L^\Gamma(\lie n^+)$ must be of the form 
\[ (x_i^-(\bar\delta)\otimes t^{-\delta})\prod_{p=1}^{q} (h_i(\bar\delta_{s})\otimes t^{-\delta_{s}}), \ \ q\le k-1.\] Now applying the element $(x_i^+(\bar\epsilon_{k+1})\otimes t^{-\epsilon_{k+1}})$ and commuting to PBW order yields terms that end in $L^\Gamma(\lie n^+)$, together with a term of length $q+1\le k$.
\end{pf}
\begin{cor}\label{h-cor} Fix $i\in I_0$ and $k_1,\ldots,k_\ell\in\bz$. Then \[\prod_{s=1}^\ell(h_i(\bar{k_s})\otimes t^{-k_s})w_\lambda\in{\text{sp}}\left\{ \prod_{s=1}^r (h_i(\bar{\ell_s})\otimes t^{-\ell_s})w_\lambda\ :\ r\le \lambda(h_i)\right\}.\]
\begin{pf}
By setting $k=\ell$ in Proposition~\ref{commute}, we see that
  \[0=\prod_{s=1}^\ell(x_i^+(\bar{k_s})\otimes t^{-k_s})(x_i^-\otimes 1)^\ell w_\lambda=\prod_{s=1}^\ell(h_i(\bar{k_s})\otimes t^{-k_s})w_\lambda + H.w_\lambda,\qquad \ell\ge\lambda(h_i)+1,\] where $H$ lies in the span of elements of the form $\prod_{s=1}^r (h_i(\bar{\ell_s})\otimes t^{-\ell_s})$ with $r<\ell$. The statement of the corollary now follows by induction on $\ell$.
\end{pf}
\end{cor}

\subsection{}
It still remains to show that the images under $\tau_\lambda$ of the elements from Lemma \ref{bagg-span} form a linearly independent subset of $\mathbb A_\lambda^\Gamma.$ Now, these are
 \[\left\{\bigotimes_{i\in I_0}\prod_{s=1}^{m_i}\sym_{\lambda(h_i)}(t^{-k_{i,s}})\ :\ 0\le m_i\le \lambda(h_i),\ k_{i,s}\in\bz\right\}.\] Since the tensor product preserves linear independence, it therefore suffices to check that for each $i\in I_0$ the set of products 
\[\left\{ \prod_{s=1}^{m_i}\sym_{\lambda(h_i)}(t^{-k_{i,s}})\ :\ 0\le m_i\le \lambda(h_i),\ k_{i,s}\in\bz\right\}\] is linearly independent. A slightly more general statement can be found in \cite{CFK10}; we reproduce it and include a proof here for convenience. Recall from Section~\ref{sym-def} that for any associative algebra $B$ we have a map $\sym_N: B\to B^{\otimes N}$ mapping \[b\mapsto \sum_{\ell=0}^{N-1}1^{\otimes \ell}\otimes b\otimes 1^{\otimes N-\ell-1}.\]

\begin{lem}
Let $b_0, b_1,\ldots \in B$ form a countable ordered basis, with $b_0=1$ and $b_r\in B_+$ for $r>1$.  Then the elements \[\left\{ \prod_{s=1}^\ell \sym_N(b_{r_s})\ :\ r_s\in\bz_+, \ell\le N\right\}\] are linearly independent in $B^{\otimes N}$.
\end{lem}
\begin{pf}
The projections onto the summand $B_+^{\otimes \ell}\otimes 1^{\otimes N-\ell}$ of the elements listed are \[ \sum_{\sigma\in S_{r_\ell}} \sigma.(b_{r_1}\otimes b_{r_2}\otimes\cdots\otimes b_{r_\ell})\otimes 1^{\otimes N-\ell},\] where $S_{\ell}$ acts in the obvious way on $B^{\otimes \ell}$. Since these are clearly linearly independent by the choice of $b_r$ as basis elements, the proof is complete.
\end{pf}

\subsection{}
In this subsection, we compare the algebras $\ba_{\lambda}$ and $\ba^{\Gamma}_{\overline{\lambda}}$. For this purpose, let us examine again the symmetrizer map defined in Section~\ref{admissible}:
$$\Sigma: \Xi \longrightarrow \Xi^{\Gamma} \ \ \ , \ \ \xi \mapsto \sum_{\sigma \in \Gamma} \sigma \circ \xi \circ \sigma^{-1}.$$
If we restrict this map to the functions of weight $\lambda \in P^+$, we obtain a map
$$\Sigma: \Xi_\lambda \longrightarrow \Xi_{\overline{\lambda}}^{\Gamma}.$$
Let $V_\xi$ be the simple $\Lg$-module associated to $\xi \in \Xi_{\lambda}$. It follows from the discussion in Section~\ref{admissible} that $V_\xi$ is a simple $\Lgg$-module if and only if $\supp(\xi)$ is admissible. On the other hand, $V_\xi$ has, viewed as a $\Lgg$-module, a unique simple quotient of highest weight $\overline{\lambda}$; in fact this quotient is isomorphic to $V_{\Sigma(\xi)}$, as we shall show in Proposition~\ref{simple-prop}.

Recall the evaluation map defined for any $\xi \in \Xi$ in \cite{CFK10} or, for $\eta \in \Xi^{\Gamma}_{\overline{\lambda}}$, in Section~\ref{taudescends}:
$$\eval_{\xi}: \bu(L(\lie h)) \longrightarrow \bc\qquad \text{and}\qquad \eval_{\eta}: \bu(\Lhg) \longrightarrow \bc.$$
The following theorem gives a natural embedding between the algebras $\ba_{\lambda}$ and $\ba^{\Gamma}_{\overline{\lambda}}$, and also gives a necessary and sufficient condition on $\lambda$ for this embedding to be surjective.
\begin{thm}\label{ba-thm} 
Let $\lambda = \sum m_i \omega_i \in P^+$, then there exists a natural injective map
$$\iota: \ba^{\Gamma}_{\overline{\lambda}} \hookrightarrow \ba_{\lambda}.$$
Furthermore, $\iota$ is surjective iff for each $i \in I, \ \lambda$  satisfies the following:

\begin{enumerate}
\item If $\Gamma.i = \{i\}$, then $m_i = 0$

\item If $m_i \neq 0$, then $m_{\sigma(i)} = 0$ for all $\sigma \in \Gamma\setminus \{ 1\}$.
\end{enumerate}
\end{thm}
\begin{pf}
We have seen that $\ba^{\Gamma}_{\overline{\lambda}} \cong \mathbb{A}^{\Gamma}_{\overline{\lambda}}$ and from \cite{CFK10}, we have \[\ba_{\lambda} \cong \mathbb{A}_{\lambda}=\bigotimes_{i\in I} (A(1)^{\otimes m_i})^{S_{m_i}}.\] All of these isomorphisms are, by construction, compatible with the embedding of $\bu(\Lhg)$ into $\bu(\Lh)$. So it remains to show that 
$$\mathbb{A}^{\Gamma}_{\overline{\lambda}} \hookrightarrow \mathbb{A}_{\lambda}.$$
It is sufficient to show this for each $i \in I_0$. Recall that we have identified $I_0$ with a subset of $I$.  We proceed with two exhaustive cases:  \\
First assume that $i \in I$ such that $\Gamma.i = \{i\}$, so that $\Gamma_i = \Gamma$. Then $A(\left| \Gamma_i \right|) \subsetneq A(1)$, so we have
$$((A(\left| \Gamma_i \right|))^{\otimes m_i})^{S_{m_i}} \subsetneq (A(1)^{\otimes m_i})^{S_{m_i}},$$
if $m_i > 0$.\\
In the other case, $\Gamma_i = \{1\}$, and we set $n_i = \sum_{\sigma \in \Gamma} m_{\sigma(i)}$.  Then we have 
$$(A(1)^{\otimes n_i})^{S_{n_i}} \subseteq \bigotimes_{\sigma \in \Gamma} (A(1)^{\otimes m_{\sigma(i)}})^{S_{m_{\sigma(i)}}},$$
with equality if and only if the right hand side consists of only one non-trivial tensor factor--i.e., $m_{\sigma(i)} = 0$ for $\sigma \neq 1$, which proves the theorem.
\end{pf}

\subsection{}\label{max-simple-sec}
Let $\lambda \in P^+$ amd $\xi \in \Xi_\lambda$. We have seen in Lemma~\ref{baggmax} that we can associate to $\xi$ a maximal ideal $\bm_\xi \in \Max \ba_{\lambda}$. The simple module $\ba_{\lambda}/\bm_\xi$ will be denoted by $\bc_\xi$. Similarly, for $\chi \in \Xi_{\overline{\lambda}}^\Gamma$ we denote the simple $\ba_{\overline{\lambda}}^{\Gamma}$-module by $\bc_\chi$.\\
Using the embedding $\ba_{\overline{\lambda}}^{\Gamma} \hookrightarrow \ba_{\lambda}$ of Theorem~\ref{ba-thm}, we see that for every $\xi \in \Xi_\lambda$, $\bc_\xi$ is a simple $\ba_{\overline{\lambda}}^{\Gamma}$-module. Then we have by construction of the symmetrizer the following:
\begin{prop}\label{simple-prop}
Let $\xi \in \Xi_\lambda$.  Then $\bc_\xi \cong \bc_{\Sigma(\xi)}$ as $\ba_{\overline{\lambda}}^{\Gamma}$-modules.
\end{prop}
\begin{pf}
In \cite[Equation (5.18)]{Lau10}, the symmetrizer map was given for multiloop algebras. It was shown that for admissible $\xi \in \Xi_\lambda$, $\bc_\xi \cong \bc_{\Sigma(\xi)}$ as $\ba_{\overline{\lambda}}^{\Gamma}$-modules. If $\xi$ is not admissible, then $V_\xi$ is not simple as a $\Lgg$-module, but has a unique simple quotient. Denote this simple quotient by $V_{\xi'}$; then $\xi' \in \Xi_\lambda$ is admissible, $\bc_\xi \cong \bc_{\xi'}$ as $\ba_{\overline{\lambda}}^{\Gamma}$-modules and $\Sigma(\xi') = \Sigma(\xi)$.
\end{pf}


\section{Local Weyl modules}\label{local-comp}
In this section, we apply the twisted Weyl functor to a special class of $\Bagg$-modules, namely the simple modules. 
We have seen that these are parametrized by $\Xi_\lambda^\Gamma$.
\begin{defn} The \textit{(twisted) local Weyl module} associated to $\chi \in \Xi_\lambda^\Gamma$ is the $\Lgg$-module
$$\bw_\lambda^{\Gamma} \bc_{\chi} \coloneqq W^{\Gamma}(\lambda) \otimes_{\Bagg} \bc_{\chi}$$
\end{defn}

One compelling reason to study the local Weyl modules is the fact that they admit the following universal property:
\begin{prop}
Let $V \in \Ob\mathcal{I}_\lambda^{\Gamma}$ such that $V$ is generated by a highest weight vector $v$ of weight $\lambda$, and suppose $\dim V_\lambda = 1$. Then there exists $\chi \in \Xi_\lambda^\Gamma$ such that the assigment $w_\lambda  \otimes 1 \mapsto v$ extends to a surjective map
$$\bw_\lambda^{\Gamma} \bc_{\chi} \twoheadrightarrow V.$$
\end{prop}
\begin{pf}
By Lemma~\ref{hwmodules} and since $V$ is generated by $v$, the assigment $w_\lambda \mapsto v$ extends to a surjective map
$$W^\Gamma(\lambda) \twoheadrightarrow V.$$
Furthermore, $V_\lambda$ is an $\Bagg$-module and since $\dim V_\lambda = 1$, this module is simple. Hence by the discussion in Section~\ref{max-simple-sec}, there exists $\chi \in \Xi_\lambda^{\Gamma}$, such that $V_\lambda \cong \bc_{\chi}$ as $\Bagg$-modules. We can deduce that the map induced by $w_\lambda \mapsto v$ factors through the kernel of $\ev_\chi$ and we have:\\
$w_\lambda \otimes 1 \mapsto v$ extends to a surjective map
$$\bw_\lambda^{\Gamma} \bc_{\chi} \twoheadrightarrow V,$$
and the proposition is proven.
\end{pf}

Local Weyl modules for twisted loop algebras have been defined before in \cite{CFS08}, as well as in \cite{FKKS11} with two different approaches. We will compare these definitions and show their equivalences; we begin by defining them for $\Lg$.

\subsection{} Let $\lambda \in P^+$ and $\xi \in \Xi_{\lambda}$. The local Weyl module associated to $\xi$, as defined in \cite{CFK10}, is
$$W(\xi) \coloneqq W(\lambda) \otimes_{\ba_{\lambda}} \bc_{\xi}.$$
Local Weyl modules had been defined previously in \cite{CP01}, but we will use the notation from \cite{CFK10}. It was shown in the aforementioned series of papers (\cite{CP01}, \cite{CL06}, \cite{FoL07}, \cite{Na11}, \cite{BN04}) that
$$\dim W(\xi) = \prod_{i \in I} (\dim W(\xi_i))^{m_i},$$
where $\lambda = \sum m_i \omega_i$ and $\xi_i$ is any element of $\Xi_{\omega_i}$. This implies that the dimension of $W(\xi)$ is independent of $\xi$, but depends only on $\lambda$. Furthermore, it has been shown (for instance in \cite{CFK10}), that $W(\xi)$ has $V_\xi$ as its unique simple quotient.
\subsection{}
In \cite{FKKS11}, local Weyl modules for $\Lgg$ were defined to be the restriction of the untwisted local Weyl module for $\Lg$; they are parametrized by equivariant finitely supported functions. We should mention, that in \cite{FKKS11} local Weyl modules were defined in a much more general context. Namely, a finite group $\Gamma$ acting freely on an affine scheme $X$ and $\lie g$ by automorphisms, which clearly includes the case of twisted loop algebras.

Specifically, let $\chi \in \Xi^{\Gamma}$ and let $\xi$ be a $\chi$-admissible function as in Section~\ref{admissible}. Then one defines by restriction the $\Lgg$-module
$$W^{\Gamma}(\chi) \coloneqq W(\xi).$$
Now, since $\xi$ is admissible, it follows that $W(\xi)$ is a cyclic $\bu(\Lgg)$-module (\cite[Theorem 4.5]{FKKS11}), and it was established in \cite[Proposition 3.5]{FKKS11} that the definition of $W^{\Gamma}(\chi)$ is independent of the choice of such $\xi$. Moreover, the modules $W^{\Gamma}(\chi)$ satisfy a universal property (\cite[Theorem 4.5]{FKKS11}) similar to the universal property of local Weyl modules for loop algebras (\cite{CP01}) and generalized current algebras (\cite[Theorem 1]{CFK10}).\\

\subsection{}
In \cite{CFS08}, local Weyl modules for the twisted loop algebra were defined by a generator $w$ and certain relations. They were parametrized by a set of $n$-tuples (where $n = \left| I_0 \right|$) of polynomials $\pi=(\pi_i)_{i\in I_0}$ with constant term 1, and we will denote these modules by $W^{\Gamma}(\mathbf{\pi})$. Their universal property was proven in \cite[Theorem 2]{CFS08}; we cite an abbreviated version here. 
\begin{thm}
Let $\lambda \in P_0^+$ and suppose that $V$ is a finite-dimensional $\Lgg$-module generated by a one-dimensional highest weight space of weight $\lambda$. Choose a vector $v_\lambda\in V_\lambda$. Then there exists an $n$-tuple $(\mathbf{\pi}_i)_{i\in I_0}$ of polynomials such that the assigment $w \mapsto v_\lambda$ extends to a surjective map of $\Lgg$-modules 
$$W^{\Gamma}(\pi) \twoheadrightarrow V.$$
\end{thm}
The following is immediate from the universal properties established above.
\begin{cor} For each $\chi \in \Xi_\lambda^{\Gamma}$, there exists a $n$-tuple of polynomials $(\mathbf{\pi})$ such that
$$\bw_{\lambda}^{\Gamma} \bc_{\chi} \cong W^{\Gamma}(\chi) \cong W^{\Gamma}(\mathbf{\pi}),$$
and vice versa.
\end{cor}

\subsection{}
In \cite{CFS08}, the dimension and character of local Weyl modules have been computed. We recall this result (\cite[Theorem 2]{CFS08}) here since it will be useful in the proof of Theorem~\ref{maintheorem}.
\begin{thm}\label{cfs-theorem} Let $\lambda \in P^+$, and $\chi \in \Xi_{\overline{\lambda}}^\Gamma$, then
\item
$$\dim W^{\Gamma}(\overline{\lambda}) \otimes_{\ba_{\overline{\lambda}}^{\Gamma}} \bc_{\chi} = \rank_{\ba_{\lambda}} W( \lambda) = \prod (\rank_{\ba_{\omega_i}} W(\omega_i))^{m_i}.$$
In particular, the dimension is independent of $\chi$ and depends only on $\overline{\lambda}$.
Moreover, the $\lie g_0$ character is also independent of $\chi$.
\end{thm}

\subsection{}
Using the fact that $\Bagg$ is a Laurent polynomial ring (Section~\ref{section-bagg}) and the fact that 
$$\dim W^{\Gamma}(\lambda) \otimes_{\Bagg} \bc_{\chi}$$
is independent of $\chi$, we can conclude a result which was previously known for untwisted loop and current algebras:
\begin{thm}\label{rank-thm}
For $\lambda \in P_0^+,\ W^{\Gamma}(\lambda)$ is a free right $\Bagg$-module with 
$$\rank_{\Bagg} W^{\Gamma}(\lambda) = \dim \bw_\lambda^{\Gamma} \bc_{\chi}$$
for some, and hence for any, $\chi \in \Xi_\lambda^\Gamma$.
\end{thm}


\section{Proof of main theorem}\label{section8}
It remains to prove the main theorem. 
\begin{thm} For $\lambda \in P_0^+$, we have
$$W^{\Gamma}(\lambda) \hookrightarrow \bigoplus_{\overline{\mu} = \lambda}  W(\mu),$$
where the map is induced by
$$w_{\lambda} \mapsto w \coloneqq \sum_{\overline{\mu} = \lambda} w_\mu.$$
\end{thm}

By construction, we have a surjective map
$$W^{\Gamma}(\lambda) \twoheadrightarrow  \bu(\Lgg).w$$
The idea of the proof is to show that both sides are free $\Bagg$-modules of the same rank. Together with the surjectivity of the above map, this will complete the proof.

\subsection{}
We have seen in Theorem~\ref{ba-thm} that $\Bagg \subset \ba_{\mu}$ is a subalgebra, for any $\mu$ satisfying $\overline{\mu} = \lambda$. It follows that $W(\mu)$ is a right module for $\Bagg$ and hence $\bigoplus_{\overline{\mu} = \lambda} W(\mu)$ is a right module for $\Bagg$. Finally, the submodule $\bu(\Lgg).w \subset \bigoplus_{\overline{\mu} = \lambda} W(\mu) $ is a right module for $\Bagg$, being a quotient of $W^\Gamma(\lambda)$.\\ 
We want to show that $\bu(\Lgg).w$ is a free $\Bagg$-module of the same rank as $W^\Gamma(\lambda)$. Now because $\Bagg$ is a polynomial algebra, in order to prove the freeness of $\bu(\Lgg).w$ it suffices to show that the dimension of
$$\bu(\Lgg).w \otimes_{\Bagg} \bc_{\chi}$$
is independent of the maximal ideal $\chi \in \Xi^{\Gamma}_{\lambda}$.\\
In order to prove this, we will need the following lemma:
\begin{lem}\label{crucial-lem}
For each $\chi \in \Xi^{\Gamma}_{\lambda}$, there exists $\tau\in P^+$ and $\xi \in \Xi_\tau$ such that $\xi$ is $\chi$-admissible and 
$$\bu(\Lgg).w \otimes_{\Bagg} \bc_{\chi} \cong_{\Lgg}  W(\tau) \otimes_{\ba_{\tau}} \bc_{\xi}.$$
\end{lem}
Assuming the lemma, we prove Theorem~\ref{maintheorem} as follows.
Observe that the dimension of the right hand side in Lemma~\ref{crucial-lem} is independent of $\xi$ and depends only on $\tau$: it is equal to the rank of $W(\tau)$ as a $\ba_{\tau}$-module. By Theorem~\ref{cfs-theorem} we know that for $\tau = \sum m_i \omega_i$, 
$$\rank_{\ba_{\tau}} W(\tau) = \prod_{i\in I} (\rank_{\ba_{\omega_i}} W(\omega_i))^{m_i}.$$
On the other hand
$$\rank_{\ba_{\omega_i}} W(\omega_i) = \rank_{\ba_{\omega_{\sigma(i)}}} W(\omega_{\sigma(i)}).$$
To see this, one may recall, that $\sigma$ is an automorphism of $\Lg$, and $W(\omega_{\sigma(i)})$ is isomorphic to the pullback of the module $W(\omega_i)$ by the automorphism $\sigma^{-1}$.\\
Using this and the rank formula for the global Weyl module, we obtain that for all $\tau_1, \tau_2 \in P^+$ with $\overline{\tau_1} = \overline{\tau_2}$, 
$$\rank_{\ba_{\tau_1}} W(\tau_1) = \rank_{\ba_{\tau_2}} W(\tau_2).$$
It follows that the dimension of $\bu(\Lgg).w \otimes_{\Bagg} \bc_{\chi}$ is independent of $\chi$, and hence  $\bu(\Lgg).w \subset \bigoplus_{\overline{\mu} = \lambda} W(\mu)$ is a projective $\Bagg$-module. Since $\Bagg$ is a polynomial ring, it now follows from the famous result of Quillen that $\bu(\Lgg).w$ is a \textit{free} $\Bagg$-module.

Together with Theorem~\ref{rank-thm}, this gives for $\overline{\tau} = \lambda$,
$$\rank_{\ba_{\tau}} W(\tau) = \rank_{\Bagg} W^{\Gamma}(\lambda).$$
We therefore conclude that the rank of $\bu(\Lgg).w \subset \bigoplus_{\overline{\mu} = \lambda} W(\mu)$ as a $\Bagg$-module is equal to the rank of  $W^{\Gamma}(\lambda)$ as a $\Bagg$-module. Since we already have a surjective map 
$$W^{\Gamma}(\lambda) \twoheadrightarrow \bu(\Lgg).w $$
and both modules are free $\Bagg$-modules, the map is an isomorphism and the theorem is proven.

\subsection{}
It remains to prove Lemma~\ref{crucial-lem}:
\begin{pf}
We start by defining projection maps $\pi_\tau$, for  $\overline{\tau} = \lambda$, onto the $\tau$-th component of $\bigoplus_{\overline{\mu} = \lambda} W(\mu)$.
$$\pi_\tau : \bigoplus_{\overline{\mu} = \lambda} W(\mu) \twoheadrightarrow W(\tau),$$
and by restriction we obtain maps
$$\pi_\tau : \bu(\Lgg).w  \longrightarrow W(\tau),$$
where $w = \sum_{\overline{\mu} = \lambda} w_\mu$. By construction, we have 
$$\pi_\tau(\bu(\Lgg).w) = \bu(\Lgg).w_\tau \subset W(\tau),$$ 
the $\Lgg$-submodule of $W(\tau)$ generated through the highest weight vector $w_\tau$.

For $\chi \in \Xi_\lambda^\Gamma$, let $\xi \in \Xi$ be a $\chi$-admissible function (whose existence is assured by the discussion in Section~\ref{admissible}) and let $\tau = \wt(\xi)$. Consider the local $\Lg$-Weyl module
$$W(\tau) \otimes_{\ba_{\tau}} \bc_\xi.$$
We see, since the support of $\xi$ is admissible, that this is a cyclic $\Lgg$-module, generated by $w_\tau \otimes \bc_\xi$. In fact, $W(\tau) \otimes_{\ba_{\tau}} \bc_\xi$ is by restriction a local Weyl module for $\Lgg$, but by construction $\bc_\chi \cong \bc_\xi$ as $\Bagg$-modules. Therefore, we have
\begin{equation}\label{eqna}
W(\tau) \otimes_{\ba_{\tau}} \bc_\xi \cong_{\Lgg} W^{\Gamma}(\lambda) \otimes_{\Bagg} \bc_\chi.
\end{equation}
Since $\bc_\chi \cong \bc_\xi$, we have the trivial isomorphism
$$\bu(\Lgg).w \otimes_{\Bagg} \bc_{\chi} \cong_{\Lgg} \bu(\Lgg).w \otimes_{\Bagg} \bc_{\xi}.$$
We use the projection map $\pi_\tau$ to obtain
$$\bu(\Lgg).w \otimes_{\Bagg} \bc_{\xi}. \twoheadrightarrow  \left( \bu(\Lgg).w_\tau \right) \otimes_{\Bagg} \bc_{\xi}.$$
Combining this projection map with the fact that $\Bagg \subset \ba_\tau$ is a subalgebra, we have as $\Lgg$-modules
\begin{equation}\label{eqnb}
\bu(\Lgg).w \otimes_{\Bagg} \bc_{\chi} \twoheadrightarrow \left( \bu(\Lgg).w_\tau \right) \otimes_{\ba_\tau} \bc_{\xi}.
\end{equation}
With the considerations above (we use that the support of $\xi$ is admissible), we obtain, that
\begin{equation}\label{eqnc}
W(\tau) \otimes_{\ba_\tau} \bc_\xi  = \left( \bu(\Lg).w_\tau \right) \otimes_{\ba_\tau} \bc_{\xi} = \left( \bu(\Lgg).w_\tau \right) \otimes_{\ba_\tau} \bc_{\xi}.
\end{equation}
Combining \ref{eqnb} and \ref{eqnc}, we obtain
$$\dim  \bu(\Lgg).w \otimes_{\Bagg} \bc_{\chi} \geq \dim W(\tau) \otimes_{\ba_\tau} \bc_\xi.$$

On the other hand, since $\bu(\Lgg).w \otimes_{\Bagg} \bc_{\chi}$ is a cyclic $\Lgg$-module, generated by the highest weight space, we have
\begin{equation}\label{eqnd}
W^{\Gamma}(\lambda) \otimes_{\Bagg} \bc_\chi \twoheadrightarrow \bu(\Lgg).w \otimes_{\Bagg} \bc_{\chi}.
\end{equation}

So we obtain 
$$\dim W^{\Gamma}(\lambda) \otimes_{\Bagg} \bc_\chi  \geq \dim  \bu(\Lgg).w \otimes_{\Bagg} \bc_{\chi}.$$
Concluding we have 
$$\dim W^{\Gamma}(\lambda) \otimes_{\Bagg} \bc_\chi  \geq \dim  \bu(\Lgg).w \otimes_{\Bagg} \bc_{\chi} \geq \dim W(\tau) \otimes_{\ba_\tau} \bc_\xi,$$
With \ref{eqna}, we conclude that we have equality throughout. That is, 
$$W^{\Gamma}(\lambda) \otimes_{\Bagg} \bc_\chi  \cong_{\Lgg} \bu(\Lgg).w \otimes_{\Bagg} \bc_{\chi} \cong_{\Lgg} W(\tau) \otimes_{\ba_\tau} \bc_\xi,$$
which finishes the proof.
\end{pf}

\subsection{}
We shall give some concluding remarks. 
\begin{rem} For $\lambda \in P_0^+$, it is not true in general that $W^\Gamma(\lambda)$ embeds into $W(\mu)$ with $\mu$ satisfying $\bar{\mu} = \lambda$. The canonical map, sending $w_\lambda$ to  $w_\mu$ is not injective.
\end{rem}
In some cases, however, this canonical map is an embedding.
\begin{rem} It can be shown that for $\omega_i \in P_0^+$, we have $W^\Gamma(\omega_i) \hookrightarrow W(\omega_j)$, for all $j \in \Gamma.i$.
\end{rem}

As mentioned in the introduction, global Weyl modules were defined only for $\lie g \otimes A$. In this paper, we have now taken a first step toward equivariant map algebras. In the twisted loop case, we still have a Cartan subalgebra, so we have weights and can define the global Weyl modules by generators and relations. In the general situation, there might be no non-zero Cartan subalgebra, so the existence of weights is no longer assured. The lack of a canonical triangular decomposition also makes it unclear how to define global Weyl modules by generators and relations. Using our result for twisted loop algebras (the global Weyl module is a submodule of an appropriate direct sum of global Weyl modules for $\Lg$) one may attempt to define global Weyl modules for $(\lie g \otimes A)^{\Gamma}$ as submodules in a direct sum of global Weyl modules for $\lie g \otimes A$. Proving that they admit sufficient properties to justify the name \textit{global Weyl module} would the main task in such a project.


\bibliographystyle{alpha}
\bibliography{twistedglobal-biblist}

\begin{thebibliography}{CFK10}

\bibitem[BN04]{BN04}
Jonathan Beck and Hiraku Nakajima.
\newblock Crystal bases and two-sided cells of quantum affine algebras.
\newblock {\em Duke Math. J.}, 123(2):335--402, 2004.

\bibitem[Car05]{Ca05}
Roger Carter.
\newblock {\em Lie algebras of finite and affine type}, volume~96 of {\em
  Cambridge Studies in Advanced Mathematics}.
\newblock Cambridge University Press, Cambridge, 2005.

\bibitem[CFK10]{CFK10}
Vyjayanthi Chari, Ghislain Fourier, and Tanusree Khandai.
\newblock A categorical approach to {W}eyl modules.
\newblock {\em Transform. Groups}, 15(3):517--549, 2010.

\bibitem[CFS08]{CFS08}
Vyjayanthi Chari, Ghislain Fourier, and Prasad Senesi.
\newblock Weyl modules for the twisted loop algebras.
\newblock {\em J. Algebra}, 319(12):5016--5038, 2008.

\bibitem[CL06]{CL06}
Vyjayanthi Chari and Sergei Loktev.
\newblock Weyl, {D}emazure and fusion modules for the current algebra of
  {$\mathfrak{sl}_{r+1}$}.
\newblock {\em Adv. Math.}, 207(2):928--960, 2006.

\bibitem[CP01]{CP01}
Vyjayanthi Chari and Andrew Pressley.
\newblock Weyl modules for classical and quantum affine algebras.
\newblock {\em Represent. Theory}, 5:191--223 (electronic), 2001.

\bibitem[FK]{FK11}
Ghislain Fourier and Deniz Kus.
\newblock Demazure modules and {W}eyl modules: {T}he twisted current case.
\newblock arXiv:1108.5960.

\bibitem[FKKS]{FKKS11}
Ghislain Fourier, Tanusree Khandai, Deniz Kus, and Alistair Savage.
\newblock Local {W}eyl modules for equivariant map algebras with free abelian
  group actions.
\newblock arXiv:1103.5766.

\bibitem[FL04]{FL04}
Boris Feigin and Sergei Loktev.
\newblock Multi-dimensional {W}eyl modules and symmetric functions.
\newblock {\em Comm. Math. Phys.}, 251(3):427--445, 2004.

\bibitem[FL07]{FoL07}
Ghislain Fourier and Peter Littelmann.
\newblock Weyl modules, {D}emazure modules, {KR}-modules, crystals, fusion
  products and limit constructions.
\newblock {\em Adv. Math.}, 211(2):566--593, 2007.

\bibitem[Her10]{Her10}
David Hernandez.
\newblock Kirillov-{R}eshetikhin conjecture: the general case.
\newblock {\em Int. Math. Res. Not. IMRN}, (1):149--193, 2010.

\bibitem[Hum72]{Hu72}
James~E. Humphreys.
\newblock {\em Introduction to {L}ie algebras and representation theory}.
\newblock Springer-Verlag, New York, 1972.
\newblock Graduate Texts in Mathematics, Vol. 9.

\bibitem[Kac90]{K90}
Victor~G. Kac.
\newblock {\em Infinite-dimensional {L}ie algebras}.
\newblock Cambridge University Press, Cambridge, third edition, 1990.

\bibitem[Lau10]{Lau10}
Michael Lau.
\newblock Representations of multiloop algebras.
\newblock {\em Pacific J. Math.}, 245(1):167--184, 2010.

\bibitem[Nao]{Na11}
Katsuyuki Naoi.
\newblock Weyl modules, {D}emazure modules and finite crystals for non-simply
  laced type.
\newblock arXiv:1012.5480.

\bibitem[NSS]{NSS}
Erhard Neher, Alistair Savage, and Prasad Senesi.
\newblock Irreducible finite-dimensional representations of equivariant map
  algebras.
\newblock \emph{Trans. Amer. Math. Soc.} (to appear), arXiv:0906.5189.

\end{thebibliography}
\end{document}